\documentclass[12pt]{amsart}
\usepackage{lipsum}
\usepackage{amssymb,amsthm}
\usepackage{amsmath}
\usepackage[foot]{amsaddr}
\usepackage{tikz}
\usepackage{graphicx}
\usepackage{subcaption}
\usepackage{caption}
\usepackage[shortlabels]{enumitem}
\usetikzlibrary{arrows}
\usepackage{float}
\usepackage{graphicx}
\usepackage{subcaption}
\usepackage{hyperref}
\usepackage{pgfplots}

\newcommand\irregularcircle[2]{% radius, irregularity
  \pgfextra {\pgfmathsetmacro\len{(#1)+rand*(#2)}}
  +(0:\len pt)
  \foreach \a in {10,20,...,350}{
    \pgfextra {\pgfmathsetmacro\len{(#1)+rand*(#2)}}
    -- +(\a:\len pt)
  } -- cycle
}

\graphicspath{{NewCompIdeasFigs/}}
\usetikzlibrary{decorations.pathmorphing}
\tikzset{snake it/.style={decorate, decoration=snake}}

\makeatletter
\newcommand*{\rom}[1]{\expandafter\@slowromancap\romannumeral #1@}
\makeatother

\numberwithin{equation}{section}

\theoremstyle{plain}
\newtheorem{theorem}{Theorem}
\numberwithin{theorem}{section}
\newtheorem{proposition}[theorem]{Proposition}

\newtheorem{corollary}[theorem]{Corollary}
\theoremstyle{definition}
\newtheorem{definition}[theorem]{Definition}
\theoremstyle{remark}
\newtheorem{remark}[theorem]{Remark}

\theoremstyle{remark}

\theoremstyle{remark}

\newcommand{\smo}{\setminus \{\mathbf{0}\}}

\newcommand{\norm}[1]{\left\lVert#1\right\rVert}      % Norm
\newcommand{\abs}[1]{\left|#1\right|}                 % Absolutbetrag
\newcommand{\paren}[1]{\left(#1\right)}               % Klammern
\newcommand{\bparen}[1]{\left[#1\right]}               % eckige Klammern
\newcommand{\sparen}[1]{\left\{#1\right\}}      % Mengenklammer
\renewcommand{\d}{\,\mathrm{d}}  % within the integral sign

\newcommand{\dd}{\mathrm{d}}  % without the space
\newcommand{\supp}{\operatorname{supp}} % singular support

\newcommand{\Ac}{\mathcal{A}}

\newcommand{\Cc}{\mathcal{C}}

\newcommand{\Dc}{\mathcal{D}}
\newcommand{\Ec}{\mathcal{E}}
\newcommand{\Fc}{\mathcal{F}}

\newcommand{\Lc}{\mathcal{L}}

\newcommand{\Rc}{\mathcal{R}}

\newcommand{\Sc}{\mathcal{S}}
\newcommand{\Tc}{\mathcal{T}}

\newcommand{\WF}{\mathrm{WF}}                         % Wavefront set
\newcommand{\wf}{\mathrm{WF}}                         % Wavefront set

\newcommand{\ve}{\mathbf{e}}
\newcommand{\inv}{^{-1}}

\newcommand{\partyf}[2]{\frac{\partial #2}{\partial y_{#1}}}

\newcommand{\bpm}{\begin{pmatrix}}
\newcommand{\epm}{\end{pmatrix}}

\newcommand{\xo}{x_0}
\newcommand{\yo}{y_0}
\newcommand{\xoyo}{(x_0,y_0)}

\newcommand{\vx}{{\mathbf{x}}}
\newcommand{\vxo}{{{\mathbf{x}_0}}}

\newcommand{\veth}{\mathbf{e}_3}

\newcommand{\vy}{{\mathbf{y}}}
\newcommand{\vz}{{\mathbf{z}}}
\newcommand{\vs}{\mathbf{s}}

\newcommand{\vxi}{{\boldsymbol{\xi}}}
\newcommand{\vxio}{{\boldsymbol{\xi}_0}}

\newcommand{\vsig}{{\boldsymbol{\sigma}}} %when $\sigma$ is one dimensional, 
\newcommand{\oB}{\overline{B}}

\newcommand{\tY}{\widetilde{Y}}

\newcommand{\otp}{[0,2\pi]}

\newcommand{\opit}{0,\pi/2} %note here there is no parenthesis around
\newcommand{\tr}{\operatorname{tr}}

\newcommand{\rpt}{\mathbb{RP}^2}

\newcommand{\ab}{(\alpha,\beta)}

\newcommand{\rr}{{{\mathbb R}}}

\newcommand{\rtwo}{{{\mathbb R}^2}}
\newcommand{\rthree}{{{\mathbb R}^3}}
\newcommand{\rn}{{{\mathbb R}^n}}

\newcommand{\drn}{{\dot{{\mathbb R}^n}}}

\newcommand{\st}{\hskip 0.3mm : \hskip 0.3mm}

\newcommand{\be}{\begin{equation}}
\newcommand{\ee}{\end{equation}}
\newcommand{\bea}{\begin{eqnarray}}
\newcommand{\eea}{\end{eqnarray}}
\newcommand{\bean}{\begin{eqnarray*}}
\newcommand{\eean}{\end{eqnarray*}}

\newcommand{\bel}[1]{\begin{equation}\label{#1}}

\newcommand{\eel}[1]{{\label{#1}\end{equation}}}

\usepackage{fullpage}

\usepackage{kbordermatrix}
\renewcommand{\kbldelim}{(}% Left delimiter
\renewcommand{\kbrdelim}{)}% Right delimiter

\usepackage{datetime}
\usetikzlibrary{quotes, angles}

\title[short]{Microlocal properties of seven-dim\-en\-sion\-al lemon and apple Radon transforms with applications in Compton scattering tomography\\{\footnotesize\ddmmyyyydate\today~\currenttime}}
\author{James W. Webber\textsuperscript{$\dagger$}}
\author{Eric Todd Quinto\textsuperscript{$\ddagger$}}
\address[James W. Webber (corresponding author)]{Department of Obstetrics and Gynecology, Brigham and Women's Hospital, 221 Longwood Ave. Boston, MA 02115}
\address[Eric Todd Quinto]{Department
of Mathematics, Tufts University, Medford, MA USA\\Partial support from NSF grant DMS
1712207 and Simons Foundation  grant 70855 }
\email[A1,A2]{jwebber5@bwh.harvard.edu\textsuperscript{$\dagger$} and Todd.Quinto@tufts.edu\textsuperscript{$\ddagger$}}
\begin{document}

\begin{abstract}
We present a microlocal analysis of two novel Radon transforms of interest in Compton Scattering Tomography (CST),
which map compactly supported
$L^2$ functions to their integrals over seven-dimensional sets of apple and lemon
surfaces. Specifically, we show that the apple and  lemon transforms
are elliptic Fourier Integral Operators (FIO), which satisfy the
Bolker condition. After an analysis of the full seven-dimensional
case, we focus our attention on $n$-D subsets of apple and lemon surfaces
with fixed central axis, where $n<7$. Such subsets of surface integrals have
applications in airport baggage and security screening. When the data dimensionality is restricted, the
apple transform is shown to violate the Bolker condition, and there
are artifacts which occur on apple-cylinder intersections. The lemon transform is shown to satisfy the Bolker
condition, when the support of the function is restricted to the strip $\sparen{0<z<1}$.
%upper half space.
\end{abstract}

\maketitle

\section{Introduction}\label{sect:intro} 
\begin{figure}
\centering
\begin{tikzpicture}[scale=2]
\draw   [blue] (0.7,0) circle (1);
\draw   [blue] (-0.7,0) circle (1);
\draw  [dashed]  (0,0) circle (0.7141);
\draw [thick, ->] (-2,0)--(2,0)node[below] {$x$};
\draw [thick, ->] (0,-2)--(0,2)node[left] {$y$};
\draw [->] (0.75,0.5)node[right] {lemon cross-section}--(0.3,0);
\draw [->] (-2,0.5)node[left] {apple cross-section}--(-1.7,0);
\end{tikzpicture}
\caption{2-D cross section of a spindle torus centered on the origin,
with axis of revolution $y$. The lemon cross-section is the intersection of the interior of the dashed circle with the torus cross-section. The apple cross section is the intersection of the torus cross-section with the exterior of the dashed circle. The lemon/apple is the surface of revolution of the lemon/apple cross-section about $y$.}
\label{fig0}
\end{figure}
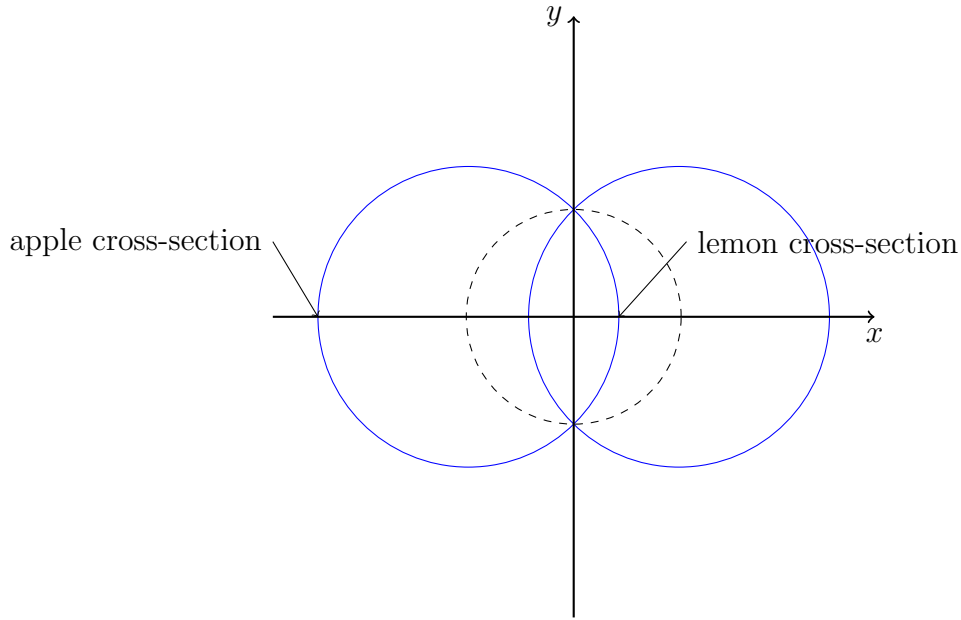

In this paper, we present a
novel microlocal analysis of two Radon transforms of interest in
CST, which take the integrals of a
function over seven dimensional sets of lemon and apple surfaces. A
``lemon" (also called a ``spindle" in some works
\cite{me2,rigaud20183d,webberholman}) refers to the interior part of a
spindle (or self-intersecting) torus, and an ``apple" is the exterior.
See figure \ref{fig0} for a 2-D cross-section of a spindle torus, where
we have highlighted the lemon and apple parts. 
The literature considers lemon and apple transforms in 3-D CST \cite{webberholman, me2, rigaud20183d, webber2019compton, rigaud20213d, rigaud2021reconstruction, cebeiro2021three}, where the goal is to reconstruct an electron density map from Compton scattered photons. There is also a growing interest in the literature in Emission CST (ECST) \cite{3D1,3D2,nguyen2005radon,moon2017analytic,kuchment2016three}, where the aim is to reconstruct a gamma ray source from cone integral data.

In \cite{rigaud2021reconstruction}, two fixed-source CST configurations, with spherical and cylindrical detector arrays, are considered. In both cases, the data is three dimensional, and consists of a two-dimensional detector coordinate and a one-dimensional energy variable. Due to limited energy resolution, the fixed source position, and the shape of the detector surface, the data is incomplete. For example, the cylindrical acquisition geometry suffers limited angle issues. In such cases of limited data, the reconstruction becomes unstable, and there are image artifacts. The authors go on to develop a modified Kaczmarz algorithm to combat the reconstruction artifacts and test their algorithm on simulated examples with Poisson noise. Similar reconstruction instabilities can be seen also in, e.g., conventional X-ray CT with limited angle data \cite{borg2018analyzing,Ks1997:limited-angle,krishnan2014microlocal}.

\begin{figure}
\centering
\begin{tikzpicture}[scale=2.5]
\draw  [dashed] (0,0) circle (1);
\draw [thick, ->] (-2,0)--(2,0)node[below] {$x$};
\draw [thick, ->] (0,-2)--(0,2)node[left] {$y$};
\draw [blue] (-2*0.95,-2*0.95) circle (2*0.5);
\draw [blue] (-2*0.4,-2*0.4) circle (2*0.5);
\fill[white] (-1.35,-1.35) circle (0.628);
%\node at (-1.35,-1.35) {$\vx_0$};

\draw [blue, domain=30:150] plot ({cos(\x)}, {sin(\x)+0.6});
\draw [blue, domain=-0.866:0.866] plot(\x, {-sqrt(1-pow(\x,2))+1.6});
\coordinate (c) at (0,0);
\draw[red,rounded corners=1mm] (c) \irregularcircle{0.75cm}{1mm};

\node at (0.3,0.3) {$f$};
\draw [->] (1.1,-1.1)node[below] {open unit ball}--(0.7,-0.7);
\draw [thick, ->] (-1.35,-1.35)node[below] {$\vx^{(2)}_0$}--(-2.025,-0.675)node[left] {$\vxi_2$};
\draw [thick, ->] (0,1.1)node[left] {$\vx^{(1)}_0$}--(1.2,1.1)node[above] {$\vxi_1$};
\draw [<->] (-1.35,-1.35)--(-0.8,-0.8);
\node at (-1.1,-1) {$t$};
\draw [<->] (-1.507,-0.0929)--(-0.8,-0.8);
\node at (-1.25,-0.5) {$r$};

\draw [->] (-0.6,-2.4)node[right] {apple cross-section}--(-1.19,-2.607);
\draw [->] (-1.4,0.8)node[left] {lemon cross-section}--(-0.86,1.1);
\end{tikzpicture}
\caption{}
\label{fig1}
\end{figure}

In \cite{webberholman}, the authors present a microlocal analysis of the lemon transform introduced in \cite{me2}. The acquisition geometry consists of a single rotating source and detector on a fixed axis. As in \cite{rigaud2021reconstruction}, the data is three-dimensional, and, in this case, consists of a 2-D rotation and a 1-D energy variable. The lemon transform is shown to violate the Bolker condition, and there are artifacts induced by flowout which appear as a spherical blurring effect in the reconstruction. There are also invisible singularities near the origin due to limited energy resolution. In \cite{me2}, an algebraic reconstruction method is proposed to invert the lemon transform. Here artifacts are observed in reconstructions with noisy data, in line with the theory of \cite{webberholman}.

In \cite{cebeiro2021three}, the authors introduce a scanning modality in 3-D CST using a fixed source and single rotating detector restricted to a spherical surface. The data, in this case, has three degrees of freedom, and consists of a 2-D detector rotation and a 1-D energy variable. The authors model the Compton scatter intensity using a new apple Radon transform, and they derive an explicit inversion formula using a spherical harmonic expansion and Volterra integral equation theory. Additionally, a hybrid analytic/algebraic reconstruction algorithm is presented and tested on simulated phantoms with added pseudo random noise. The authors discover blurring artifacts in the reconstructions, which indicate instabilities due to limited data, as is, for example, discovered in \cite{rigaud2021reconstruction}.

In the works discussed above, a number of imaging modalities are
introduced based on practical machine designs, and the data dimension
is such that the reconstruction target is determined. That is, the
reconstruction target and data are both three-dimensional. The set of
spindle tori in 3-D space is seven-dimensional, and hence the
literature thus far considers only limited data problems in CST, i.e.,
3-D subsets of the full 7-D set of tori are considered. This often
leads to artifacts and instabilities in the reconstruction due to, for
example, limited angles (as in \cite{rigaud2021reconstruction}) and
failure to satisfy the Bolker condition \cite{webberholman}. In this
paper, we wish to investigate the problem instability and presence of
artifacts when there are no limits to the data dimensionality in CST,
and we have knowledge of a seven-dimensional set of apple and lemon
integrals in 3-D space. This can be considered a best case scenario in
CST in terms of data dimensionality. Specifically, we consider the
scanning geometry illustrated in figure \ref{fig1}. Here, we have
shown an $(x,y)$ plane cross-section of the scanning geometry. The
scanning target ($f$) is supported on the open unit ball and is
illustrated by an uneven red boundary. Example lemon and apple cross
sections are drawn in blue, with centers $\vx^{(1)}_0$ and
$\vx^{(2)}_0$, and axis of rotation $\vxi_1$ and $\vxi_2$,
respectively. The apple radius is denoted by $r$, and the distance
from $\vx^{(2)}_0$ to the center of the apple tube is denoted by $t$.
We consider the apple and lemon surfaces whose points of
self-intersection (which we will call \emph{singular points}) lie outside the
open unit ball. We do this to avoid singularities in the apple/lemon
surface measure. In CST, the singular points of the lemons
and apples correspond to source and detector coordinates. So, in the context
of CST, our geometry consists of all source and detector positions
which lie outside the unit ball (this is a six-dimensional set).
Additionally, we can vary the torus radius ($r$), which in CST is equivalent to the photon energy \cite{rigaud2021reconstruction}. Thus, in total, our
data set is seven-dimensional.

%\tc{Edit the red text below to reflect what we actually prove about
%ellipticity, etc.}

Motivated by the geometry of figure \ref{fig1}, we introduce novel lemon and apple Radon transforms, which map $f$ to its integrals over seven-dimensional sets of apple and lemon surfaces. Our main theorem proves that the lemon and apple transforms are elliptic FIO which satisfy the Bolker condition. 
%\tred{Following this, we show that the normal operators of the lemon and apple transforms are elliptic Pseudodifferential Operators (PDO), and further we show the boundedness of solution in Sobolev space, as is done, e.g., for the classical Radon transform in \cite[Theorem 5.1]{natterer}.}
Additionally, we consider the practical applications of our theory to other scanning geometries from the literature. Specifically, we consider the scanning geometry of \cite{WebberQuinto2020II}, which is designed for use in airport baggage screening, and discuss the microlocal properties of lemon and apple transforms which induce translation on the scanning target.

The remainder of this paper is organized as follows. In section \ref{sect:defns}, we give some preliminary definitions and theorems that will be used in our analysis. In section \ref{radon section}, we introduce novel lemon and apple transforms, which map compactly supported $L^2$ functions to their integrals over seven-dimensional sets of lemon and apple surfaces, respectively, as pictured in figure \ref{fig1}. Here we prove our main theorem, which shows that the lemon and apple transforms are elliptic FIO which satisfy the Bolker condition. In section \ref{CST}, we consider a practical scanning geometry in CST, first introduced in \cite{WebberQuinto2020II}, and discuss the artifacts in lemon and apple integral reconstructions when the axis of revolution of the lemons$\backslash$apples is fixed, and the target function undergoes a 2-D translation. 
\section{Definitions and preliminary theorems}\label{sect:defns} 

We next provide some
notation and definitions.  Let $X$ and $Y$ be open subsets of
$\rn$.  Let $\Dc(X)$ be the space of smooth functions compactly
supported on $X$ with the standard topology and let $\mathcal{D}'(X)$
denote its dual space, the vector space of distributions on $X$.  Let
$\Ec(X)$ be the space of all smooth functions on $X$ with the standard
topology and let $\mathcal{E}'(X)$ denote its dual space, the vector
space of distributions with compact support contained in $X$. Finally,
let $\Sc(\rn)$ be the space of Schwartz functions, that are rapidly
decreasing at $\infty$ along with all derivatives. See \cite{Rudin:FA}
for more information.

{For a function $f$ in the Schwartz space $\Sc(\mathbb{R}^n)$ or in
$L^2(\rn)$, we use $\mathcal{F}f$ and $\mathcal{F}^{-1}f$ to denote
the Fourier transform and inverse Fourier transform of $f$,
respectively (see \cite[Definition 7.1.1]{hormanderI}).  Note that
$\Fc\inv \Fc f(\vx)= \frac{1}{(2\pi)^n}\int_{\vy\in\rn}\int_{\vz\in
\rn} \exp((\vx-\vz)\cdot \vy)\,
f(\vz)\d \vz\d \vy$.}

We use the standard multi-index notation: if
$\alpha=(\alpha_1,\alpha_2,\dots,\alpha_n)\in \sparen{0,1,2,\dots}^n$
is a multi-index and $f$ is a function on $\rn$, then
\[\partial^\alpha f=\paren{\frac{\partial}{\partial
x_1}}^{\alpha_1}\paren{\frac{\partial}{\partial
x_2}}^{\alpha_2}\cdots\paren{\frac{\partial}{\partial x_n}}^{\alpha_n}
f.\] If $f$ is a function of $(\vy,\vx,\vs)$ then $\partial^\alpha_\vy
f$ and $\partial^\alpha_\vs f$ are defined similarly.

  We identify cotangent
spaces on Euclidean spaces with the underlying Euclidean spaces, so we
identify $T^*(X)$ with $X\times \rn$.

If $\phi$ is a function of $(\vy,\vx,\vs)\in Y\times X\times \rr^N$
then we define $\dd_{\vy} \phi = \paren{\partyf{1}{\phi},
\partyf{2}{\phi}, \cdots, \partyf{n}{\phi} }$, and $\dd_\vx\phi$ and $
\dd_\vs \phi $ are defined similarly. We let $\dd\phi =
\paren{\dd_{\vy} \phi, \dd_{\vx} \phi,\dd_\vs \phi}$.

% \[\begin{gathered}\dd_{\vy} \phi = \paren{\partyf{1}{\phi},
% \partyf{2}{\phi}, \cdots, \partyf{n}{\phi} },\ \dd_\vs \phi =
% \paren{\partsif{1}{\phi},\partsif{2}{\phi}, \cdots, \partsif{N}{\phi}
% }\\ \text{ and }\ \dd\phi(\vx,\vs) = \paren{\dd_{\vy} \phi(\vy,
% \vx,\vs), \dd_{\vx} \phi(\vy,\vx,\vs),\dd_\vs
% \phi(\vy,\vx,\vs)}\in \rn\times \rn\times\rr^N.\end{gathered}\]

We use the convenient notation that if $A\subset \rr^m$, then $\dot{A}
= A\smo$.

The singularities of a function and the directions in which they occur
are described by the wavefront set \cite[page
16]{duistermaat1996fourier}: 
\begin{definition}
\label{WF} Let $X$ Let an open subset of $\rn$ and let $f$ be a
distribution in $\mathcal{D}'(X)$.  Let $(\vx_0,\vxi_0)\in X\times
\drn$.  Then $f$ is \emph{smooth at $\vx_0$ in direction $\vxio$} if
there exists a neighborhood $U$ of $\vx_0$ and $V$ of $\vxi_0$ such
that for every $\phi\in \Dc(U)$ and $N\in\mathbb{R}$ there exists a
constant $C_N$ such that for all $\vxi\in V$,
\begin{equation}
\left|\Fc(\phi f)(\lambda\vxi)\right|\leq C_N(1+\abs{\lambda})^{-N}.
\end{equation}
The pair $(\vx_0,\vxio)$ is in the \emph{wavefront set,} $\wf(f)$, if
$f$ is not smooth at $\vx_0$ in direction $\vxio$.
\end{definition}
 This definition follows the intuitive idea that the elements of
$\WF(f)$ are the point--normal vector pairs above points of $X$ at
which $f$ has singularities.  For example, if $f$ is the
characteristic function of the unit ball in $\mathbb{R}^3$, then its
wavefront set is $\WF(f)=\{(\vx,t\vx): \vx\in S^{2}, t\neq 0\}$, the
set of points on a sphere paired with the corresponding normal vectors
to the sphere.

%\begin{equation}

%\end{equation}
%That is, 

The wavefront set of a distribution on $X$ is normally defined as a
subset the cotangent bundle $T^*(X)$ so it is invariant under
diffeomorphisms, but we do not need this invariance, so we will
continue to identify $T^*(X) = X \times \rn$ and consider $\WF(f)$ as
a subset of $X\times \drn$.

%Let $X$ and $Y$ be open subsets of $\rn$, $m \in\mathbb{R}$.

 \begin{definition}[{\cite[Definition 7.8.1]{hormanderI}}] \label{ellip}We define
 $S^m(Y\times X\times \mathbb{R}^N)$ to be the
set of $a\in \Ec(Y\times X\times \mathbb{R}^N)$ such that for every
compact set $K\subset Y\times X$ and all multi--indices $\alpha,
\beta, \gamma$ the bound
\[
\left|\partial^{\gamma}_{\vy}\partial^{\beta}_{\vx}\partial^{\alpha}_{\vsig}a(\vy,\vx,\vsig)\right|\leq
C_{K,\alpha,\beta,\gamma}(1+\norm{\vsig})^{m-|\alpha|},\ \ \ (\vy,\vx)\in K,\
\vsig\in\mathbb{R}^N,
\]
holds for some constant $C_{K,\alpha,\beta,\gamma}>0$. 

 The elements of $S^m$ are called \emph{symbols} of order $m$.  Note
that these symbols are sometimes denoted $S^m_{1,0}$.  The symbol
$a\in S^m(Y,X,\rr^N)$ is \emph{elliptic} if for each compact set
$K\subset Y\times X$, there is a $C_K>0$ and $M>0$ such that
\bel{def:elliptic} \abs{a(\vy,\vx,\vsig)}\geq C_K(1+\norm{\vsig})^m,\
\ \ (\vy,\vx)\in K,\ \norm{\vsig}\geq M.\ee 
\end{definition}

\begin{definition}[{\cite[Definition
        21.2.15]{hormanderIII}}] \label{phasedef}
A function $\phi=\phi(\vy,\vx,\vsig)\in
\Ec(Y\times X\times\dot{\mathbb{R}^N})$ is a \emph{phase
function} if $\phi(\vy,\vx,\lambda\vsig)=\lambda\phi(\vy,\vx,\vsig)$, $\forall
\lambda>0$ and $\mathrm{d}\phi$ is nowhere zero. The
\emph{critical set of $\phi$} is
\[\Sigma_\phi=\{(\vy,\vx,\vsig)\in Y\times X\times\dot{\mathbb{R}^N}
: \dd_{\vsig}\phi=0\}.\] 
 A phase function is
\emph{clean} if the critical set $\Sigma_\phi = \{ (\vy,\vx,\vsig) \ : \
\mathrm{d}_\vsig \phi(\vy,\vx,\vsig) = 0 \}$ is a smooth manifold {with tangent space defined as the kernel of $d\,(d_\sigma\phi)$ on $\Sigma_\phi$. Here, the derivative $\mathrm{d}$ is applied component-wise to the vector-valued function $d_\sigma\phi$. So, $d\,(d_\sigma\phi)$ is treated as a Jacobian matrix of dimensions $N\times (2n+N)$.}
\end{definition}
\noindent By the {Constant Rank Theorem} the requirement for a phase
function to be clean is satisfied if
$\mathrm{d}\paren{\mathrm{d}_\vsig
\phi}$ has constant rank.

\begin{definition}[{\cite[Definition 21.2.15]{hormanderIII} and
      \cite[section 25.2]{hormander}}]\label{def:canon} Let $X$ and
$Y$ be open subsets of $\rn$. Let $\phi\in \Ec\paren{Y \times X \times
{\rr}^N}$ be a clean phase function.  In addition, we assume that
$\phi$ is \emph{nondegenerate} in the following sense:
\[\text{$\dd_{\vy}\phi$ and $\dd_{\vx}\phi$ are never zero on
$\Sigma_{\phi}$.}\]
  The
\emph{canonical relation parametrized by $\phi$} is defined as
\begin{equation}\label{def:Cgenl} \begin{aligned} \Cc=&\sparen{
\paren{\paren{\vy,\dd_{\vy}\phi(\vy,\vx,\vsig)};\paren{\vx,-\dd_{\vx}\phi(\vy,\vx,\vsig)}}:(\vy,\vx,\vsig)\in
\Sigma_{\phi}},
% &\hspace{1.5cm} \vs\in \rr^N\smo,   
\end{aligned}
\end{equation}
\end{definition}

\begin{definition}\label{FIOdef}
Let $X$ and $Y$ be open subsets of $\rn$. {Let an operator $A :
\Dc(X)\to \mathcal{D}'(Y)$ be defined by the distribution kernel
$K_A\in \mathcal{D}'(X\times Y)$, in the sense that
$Af(\vy)=\int_{X}K_A(\vx,\vy)f(\vx)\mathrm{d}\vx$. Then we call $K_A$
the \emph{Schwartz kernel} of $A$}. A \emph{Fourier
integral operator (FIO)} of order $m + N/2 - n/2$ is an operator
$A:\Dc(X)\to \mathcal{D}'(Y)$ with Schwartz kernel given by an
oscillatory integral of the form
\begin{equation} \label{oscint}
K_A(\vy,\vx)=\int_{\mathbb{R}^N}
e^{i\phi(\vy,\vx,\vsig)}a(\vy,\vx,\vsig) \mathrm{d}\vsig,
\end{equation}
where $\phi$ is a clean nondegenerate phase function and $a$ is a
symbol in $S^m(Y \times X \times \mathbb{R}^N)$. The \emph{canonical
relation of $A$} is the canonical relation of $\phi$ defined in
\eqref{def:Cgenl}.

The FIO $A$ is \emph{elliptic} if its symbol is elliptic.
\end{definition}

This is a simplified version of the definition of FIO in \cite[section
2.4]{duist} or \cite[section 25.2]{hormander} that is suitable for our
purposes since our phase functions are global. Because we assume phase
functions are nondegenerate, our FIO can be defined as maps from
$\Ec'(X)$ to $\Dc'(Y)$ and sometimes on larger domains. For general
information about FIOs, see \cite{duist, hormander, hormanderIII}. For
information about the Schwartz Kernel, see \cite[Theorem
5.1.9]{hormanderI}. 

.

{Let $X$ and $Y$ be
sets and let $\Omega_1\subset X$ and $\Omega_2\subset Y\times X$. The composition $\Omega_2\circ \Omega_1$ and transpose $\Omega_2^t$ of $\Omega_2$ are defined
\[\begin{aligned}\Omega_2\circ \Omega_1 &= \sparen{y\in Y\st \exists x\in \Omega_1,\
(y,x)\in \Omega_2}\\
\Omega_2^t &= \sparen{(x,y)\st (y,x)\in \Omega_2}.\end{aligned}\]}

The H\"ormander-Sato Lemma  provides the relationship between the
wavefront set of distributions and their images under FIO.

\begin{theorem}[{\cite[Theorem 8.2.13]{hormanderI}}]\label{thm:HS} Let $f\in \Ec'(X)$ and
let ${A}:\Ec'(X)\to \Dc'(Y)$ be an FIO with canonical relation $\Cc$.
Then, $\wf({A}f)\subset \Cc\circ \wf(f)$.\end{theorem}

\begin{definition}
\label{defproj} Let $\Cc\subset T^*(Y\times X)$ be the canonical
relation associated to the FIO ${A}:\mathcal{E}'(X)\to
\mathcal{D}'(Y)$. We let $\Pi_L$ and $\Pi_R$ denote the natural left-
and right-projections of $\Cc$, projecting onto the appropriate
coordinates: $\Pi_L:\Cc\to T^*(Y)$ and $\Pi_R : \Cc\to T^*(X)$.
\end{definition}

Because $\phi$ is nondegenerate, the projections do not map to the
zero section.  
% 
% We have the following result from \cite{hormander}.
% \begin{proposition}
% \label{prop1}
% Let $\dim(X)=\dim(Y)$. Then at any point in $\Cc$:
% \begin{enumerate}[(i)]
% \item if one of $\Pi_L$ or $\Pi_R$ is a local diffeomorphism, then the
% other map is a local diffeomorphism (so $\Cc$ is a local canonical
% graph); 
% 
% \item if one of the projections $\Pi_R$ or $\Pi_L$ is singular at a
% point in $\Cc$, then so is the other. The type of the singularity may
% be different but both projections drop rank on the same set
% \begin{equation}
% \Sigma=\{(\vy,\eta; \vx,\vsig)\in \Cc :
% \det(\mathrm{d}\Pi_L)=0\}=\{(\vy,\eta; \vx,\vsig)\in \Cc : \det
% (\mathrm{d}\Pi_R)=0\}.
% \end{equation}
% \end{enumerate}
% \end{proposition}

{Let $A$ be an FIO with adjoint $A^*$. If $A$ satisfies our next definition, then $A^* A$ (or, if
$A$ does not map to $\Ec'(Y)$, then $A^* \psi A$ for an
appropriate cutoff $\psi$) is a pseudodifferential operator
\cite{GS1977, quinto}.}

\begin{definition}\label{def:bolker} Let
${A}:\Ec'(X)\to \Dc'(Y)$ be a FIO with canonical relation $\Cc$ then
{$A$} (or $\Cc$) satisfies the \emph{semi-global Bolker Condition} if
the natural projection $\Pi_L:\Cc\to T^*(Y)$ is an embedding
(injective immersion).\end{definition}

% \tc{I added a reference to the SDT since I had never heard of it.  This is the original article in which it was stated, but not proved.  Here is the URL that gave it to me:  \url{http://www.scientificlib.com/en/Mathematics/LX/SylvestersDeterminantTheorem.html}}

\begin{theorem}[Sylvester's Determinant Theorem (SDT) {\cite{SDT-web,SDT}}]
Let $A$ be an $m\times n$ matrix, and $B$ an $n\times m$ matrix. Then
$$\text{det}\paren{I_{m\times m}+AB}=\text{det}\paren{I_{n\times n}+BA}.$$
\end{theorem}

\section{Analysis of seven-dimensional lemon and apple Radon
transforms} \label{radon section} In this section, we present a
microlocal analysis of two new Radon transforms which map compactly supported $L^2$ functions to their
integrals over seven-dimensional sets of lemon and apple
surfaces. First, we give the defining equations for the apple
and lemon surfaces.

Spindle tori are described by their center, $\vxo\in \rthree$, their
axis of revolution, and parameters $s$ and $t$; $\sqrt{s}$ is the radius
and $t$ is the tube radius of the spindle torus. If $\ell$ is a line
through the origin parallel to the axis of revolution of a spindle torus,
then for some $\omega\in S^2$, one can write \[\ell = \rr\omega
:=\sparen{\nu\omega\st \nu\in \rr}\] and the axis of revolution of the
torus is $\vxo + \ell$. We will call this line $\ell$ the
\emph{directional axis} of the spindle torus (equivalently, of the
apple or lemon).

\begin{figure}[!h]
\centering
\begin{subfigure}{0.49\textwidth}
\includegraphics[ width=1\linewidth, height=2\linewidth, keepaspectratio]{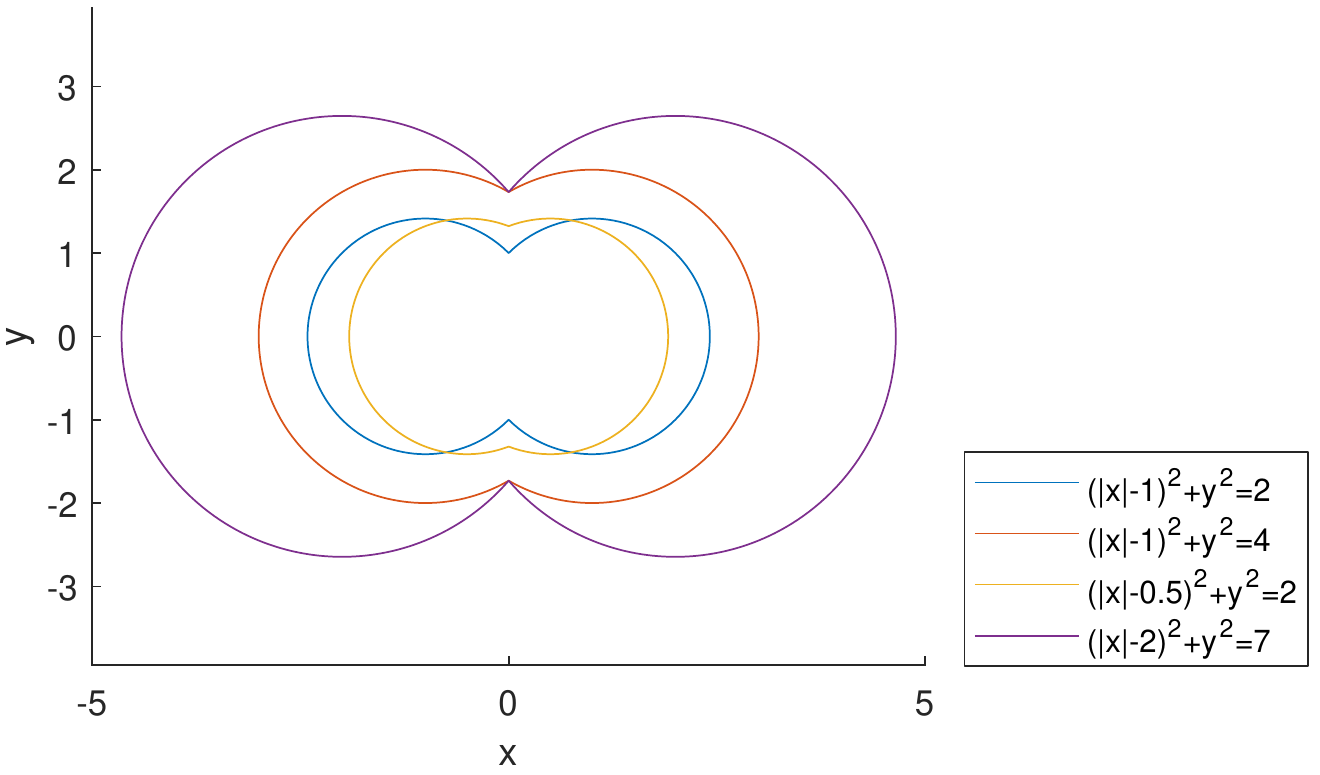} 
\subcaption{Apples ($j=1$)}
\end{subfigure}
\begin{subfigure}{0.49\textwidth}
\includegraphics[ width=1\linewidth, height=2\linewidth, keepaspectratio]{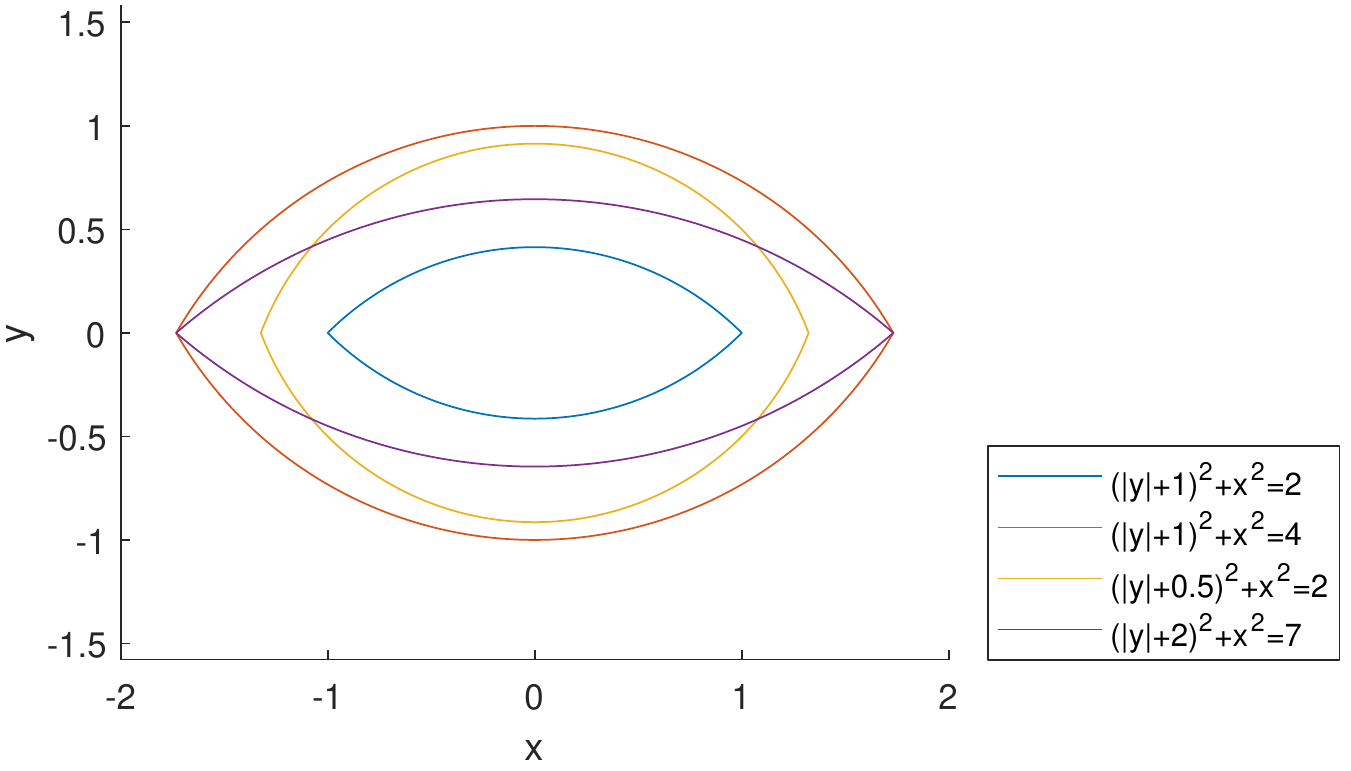} 
\subcaption{Lemons ($j=2$)}
\end{subfigure}
\caption{$(x,y)$ plane cross sections of the apple and lemon parts of a spindle
torus when $R=R\paren{0,\frac{\pi}{2}}$ (left) and $R=R\paren{\frac{\pi}{2},\frac{\pi}{2}}$ (right), $\vx_0=\textbf{0}$, and $s$ and $t$ vary between $\frac{1}{2}$ and $7$.}
\label{figlemon}
\end{figure}

We will use rotation matrices to describe the directional axes of
spindle tori. Let $(\alpha,\beta)\in \otp\times [0,\pi/2]$. Then, we
define 
\begin{equation}\label{def:R}\begin{aligned}
R=R(\alpha,\beta)& = \begin{pmatrix}
\cos\alpha &-\sin\alpha& 0 \\
\sin\alpha &\cos\alpha& 0 \\
0 &0& 1
\end{pmatrix}\begin{pmatrix}
1 &0& 0\\
0 &\cos\beta &-\sin\beta \\
0 &\sin\beta &\cos\beta 
\end{pmatrix}\\&=\begin{pmatrix}
\cos\alpha &-\sin\alpha\cos\beta& \sin\alpha\sin\beta \\
\sin\alpha &\cos\alpha\cos\beta& -\cos\alpha\sin\beta \\
0 &\sin\beta& \cos\beta
\end{pmatrix}.\end{aligned}
\end{equation} 
% \noindent\rule{\linewidth}{0.6pt}
% 
% \tc{Gulp. Here's why I'm doing this. To show Bolker, we need that $Y$
% is a manifold without boundary and this means finding an intrinsic way
% to describe all spindle torus axes.}
% 
% \tred{More details: \it Spindle tori are determined by their axis and
% then $(s,t,\vxo)$. The set of axes is the set of lines through the
% origin in $\rthree$. This set is denoted $\mathbb{RP}^2$, real
% projective space. One can also describe $\rpt$ as $S^2$ with antipodal
% points identified (so $\veth$ and $-\veth$ represent the same point in
% $\rpt$
% since they are on the same line through the origin). This would argue
% that $Y$ should be the set $(s,t,\vxo,\ell)$ where $\ell\in \rpt$. 
% 
% Then, to check Bolker, one would define local coordinates $\ab$ and
% the axis $\ell\ab$ is the axis with direction $R\ab \veth$ for $\ab\in
% \otp\times (0,\pi/2)$. 
% 
% Then, we note that we can use local coordinates by choosing a
% ``z-axis'' and then putting coordinates about that and then
% parametrizing using $R\ab$ for those axes that are not vertical or in
% the horizontal plane. Then, we choose a new ``z-axis'' and use new
% coordinates. I try to explain this at the start of the proof. 
% 
% What do you think?}
% 
% \noindent\rule{\linewidth}{0.6pt}
Let 
\bel{defns1} \vx=(x,y,z),\ \ \ \ \vx_0=(x_0,y_0,z_0),\ \ \vx_T=\vx-\vx_0,\ \ 
 \vx'=(x',y',z')=R^T(\alpha,\beta)\vx_T
 \ee
and
\bel{defns2}
h(t,\vx_0;\vx)=\norm{\vx_T}^2+t^2,\ \ \ g(\alpha,\vx_0,\beta;\vx)=\sqrt{x'^2+y'^2}.
\end{equation}
We now define
\bel{defn:Psi}\begin{aligned}
\Psi_j(s,t,\vx_0,\alpha,\beta;\vx)&:=\paren{\sqrt{x'^2+y'^2} +
(-1)^jt}^2 + z'^2-s\\
&\phantom{:}=h(t,\vx_0;\vx)+2t(-1)^jg(\vx_0,\alpha,\beta;\vx)-s.\end{aligned}
\end{equation}
and \bel{equdef}T_j(s,t,\vx_0,\alpha,\beta) = \sparen{\vx\in B\st
\Psi_j(s,t,\vx_0,\alpha,\beta;\vx)=0},\ee
for $j=1,2$. $\Psi_1$ and $\Psi_2$ are the defining equations for apple and lemon surfaces, respectively, and $T_1$ and $T_2$ are the intersections of apples and lemons with $B$. See figure \ref{figlemon} for example 2-D cross sections of
apples and lemons with the defining equations highlighted.

%Note also that $(\alpha,\beta)\mapsto R(\alpha,\beta)$ for $(\alpha,\beta)\in \otp^2$ (where we identify $0$ and $2\pi$) is an embedding into $SO(3)$. The easiest way to see this is to note that the first column of $R(\alpha,\beta)$ smoothly gives you $\alpha$ and the third row smoothly gives you $\beta$.

\subsection{Definition of apple and lemon transforms}

Throughout this paper, we let $L^2_c(X)$ denote the set of $L^2$ functions compactly supported on $X\subset\mathbb{R}^3$. Recall that the two points of intersection of the apple
(resp.\ lemon) with its axis of revolution are called the \emph{singular points} of
the apple (resp.\ lemon). Note that the singular points are the points
of intersection of the apple and lemon with the same parameters in
$Y$, and they are singular points of both the apple and lemon. We will define the apple and lemon transforms on 
functions $f\in L^2_c(X)$, where $X=B$ is the open unit ball in
$\mathbb{R}^3$, and we will need to ensure that the singular points of the
apple or lemon do not meet the closed unit ball, $\oB$. For this
reason, we define
\begin{equation}\label{def:Y}\begin{aligned}
Y=&\Large\{(s,t,\vx_0,\alpha,\beta)\in \rtwo\times \rthree \times
[0,2\pi]\times \bparen{0,\pi/2}\\
&\quad: s>t^2, \{\vx_0\pm\sqrt{s-t^2} R(\alpha,\beta)\veth\}\cap
\oB=\emptyset\Large\},\end{aligned}\end{equation} where
$\veth=(0,0,1)^T$ is the north pole. Note that every apple (for $j=1$)
and lemon (for $j=2$) with singular points not meeting $\oB$ can be
written $T_j(s,t,\vxo,\alpha,\beta)$ for some
$(s,t,\vxo,\alpha,\beta)\in Y$ because all directional axes are
generated by the map \bel{pi2 map]}\otp\times [0,\pi/2]\ni \ab\mapsto
\rr R\ab\veth. \ee

 \begin{remark}\label{rem:coords1} This map \eqref{pi2 map]} from
$\ab$ to directional axes is not injective for $\beta = 0, \pi/2$.
Therefore, we cannot use $\otp\times [\opit]$ to parameterize
direction axes, as it would cause issue later in the proofs of our main theorems.  Furthermore, our parameter space $Y$ in \eqref{def:Y} is
not a manifold without boundary because $[\opit]$ is not a manifold
without boundary. Note that we are identifying $0$ and $2\pi$ to transform $\otp$ to the manifold $S^1$.

At the start of the proof of Theorem \ref{AFIO-1}, we will define a
parameter set for spindle tori that is a manifold without boundary for
which the map to spindle tori (with singular points outside $\oB$) is
bijective. These properties are required to use the standard
microlocal analysis of Radon transforms (e.g., see \cite{GS1977}). However, we will parameterize spindle tori using $Y$ when appropriate.
\end{remark}

We define the Radon transforms which take the integrals of $f$ over
apple ($j=1$) and lemon ($j=2$) surfaces
\begin{equation}
\label{radon}
\begin{split}
\mathcal{R}_jf(s,t,\vx_0,\alpha,\beta)&=\int_{X}\norm{\nabla_{\vx}\Psi_j}\delta\paren{\Psi_j(s,t,\vx_0,\alpha,\beta;\vx)}f(\vx)\mathrm{d}\vx\\
&=\int_{-\infty}^{\infty}\int_{X}\norm{\nabla_{\vx}\Psi_j}e^{\sigma\Psi_j(s,t,\vx_0,\alpha,\beta;\vx)}f(\vx)\mathrm{d}\vx\mathrm{d}\sigma,
\end{split}
\end{equation}
and we let $$\mathcal{A}f=\mathcal{R}_1f,\ \ \ \ \ \
\mathcal{L}f=\mathcal{R}_2f$$ where $\mathcal{A}$ is called the
\emph{apple transform,} and $\mathcal{L}$ is the \emph{lemon
transform}. 

Here, we will assume the gradient of a scalar valued function is a
column vector, as are elements of $\rn$.

To ensure that $T_j(s,t,\vxo,\alpha,\beta)$
is a smooth manifold
and that the weight $\norm{\nabla_{\vx}\Psi_j}$ in
\eqref{radon} is defined, we have defined $Y$ so that it includes only the apple
and lemon surfaces whose singular points do not intersect $X=B$. This
way, in the integrals of \eqref{radon}, we stay away from the singular
points of the apples and lemons, and any singularities in the FIO
amplitudes and phases. Strictly speaking, one would add a smooth cutoff to the
symbol which is zero close to the central axis of the spindle tori, as in \cite[Lemma
3.3]{webber2021microlocal}, so the amplitude is smooth everywhere
and the phase is smooth near the support of the amplitude. However,
we do not go into such technicalities here. 

%\tred{Note that our functions have compact support away from the
%singular points on the spindle tori parameterized by $Y$ by
%definition, \eqref{def:Y}.} \medskip

Now that the apple and lemon transforms are defined we
present a separate microlocal analysis of each transform in the
following sections.

\subsection{Microlocal properties of $\mathcal{A}$; the $j=1$ case}
Here we discuss the microlocal properties of the apple transform $\mathcal{A}$. Our first theorem proves that $\mathcal{A}$ is an elliptic FIO.
\begin{theorem}
\label{AFIO-1} 
The apple transform $\mathcal{A}=\mathcal{R}_1$ of \eqref{radon} is an elliptic FIO order $-2$ from domain $\Ec'(B)$ to $\Dc'(Y)$.
\end{theorem}

\begin{proof}
To analyze $\Ac$ as an FIO, we need to parametrize apples using a
manifold without boundary, as discussed in Remark \ref{rem:coords1}.
However, $Y$ cannot be used, since it is not a manifold without
boundary since $[0,\pi/2]$ has boundary points $0,\pi/2$. To get around this, we first parametrize all spindle tori in a global
way as a manifold without boundary. This is required to use the theory
of Radon transforms as FIO \cite{GS1977}. To define this manifold, we parametrize spindle tori by points
$\vy=(s,t,\vxo,\ell)$, as discussed at the start of this section, where $\sqrt{s}$ is the radius, $t$ is the
tube radius of the spindle torus, $\vxo$ is its center, and $\ell$ is the directional axis. Recall that the directional axis of a spindle torus is the line through the origin in $\rthree$, which is parallel to the axis
of revolution of the torus, $\vxo+\ell$.  The set of lines through the
origin in $\rthree$ is denoted $\rpt$ and is called the  two-dimensional real projective space.

We let $\tY$ be the set of
$\vy=(s,t,\vxo,\ell)$ such that the singular points of the
spindle torus parameterized by $\vy$ do not meet $\oB$. Then, $\tY$
is a manifold without boundary that parameterizes all apples ($j=1$)
and all lemons ($j=2$) the singular points of which do not meet $\oB$
by the map 
\bel{tY map}\begin{gathered}\tY\ni (s,t,\vxo,\ell)\mapsto
T_j(s,t,\vxo,\alpha,\beta) \\ \text{when $\ab$ is chosen so that}\\
\ell = \rr R\ab \veth.\end{gathered}\ee Note that the map in \eqref{tY
map} and $\Ac f$ are well-defined on $\tY$ because the spindle torus
and its measure are the same no matter which $\ab$ one chooses that
satisfies $\ell = \rr R\ab\veth$. This is true by rotation invariance
of the spindle torus about its axis of revolution, $\vxo+\ell$ and
rotation invariance of the integral over the torus. Furthermore, every
spindle torus is described by a unique $(s,t,\vxo,\ell)$.

To get local coordinates on $\tY$, we need to specify local
coordinates on $\rpt$, since $(s,t,\vxo)$ are already coordinates. We
choose a vertical axis and let $\veth$ be the unit vector pointing in
the positive direction along that axis. Then, we let $\ve_1$ and
$\ve_2$ be orthogonal unit vectors so $(\ve_1,\ve_2,\veth)$ form a
right-hand coordinate system in $\rthree$. Now, we define the domain
of the coordinate map
\bel{def:Y'} Y' = \sparen{(s,t,\vxo,\alpha,\beta)\in Y\st \beta\in
(0,\pi/2)}.\ee
 Then, local coordinates on
$\tY$ are given by \bel{tY map)} Y'\ni
(s,t,\vxo,\alpha,\beta) \mapsto (s,t,\vxo,\rr R\ab \veth).\ee
For different choices of basis $(\ve_1,\ve_2,\ve_3)$ on $\rthree$ with vertical axis in direction of  $\ve_3$, this coordinate map describes a coordinate chart on $\tY$. 

We will work in these coordinates and use the notation
\eqref{def:R}-\eqref{defn:Psi}, and \eqref{radon} for the rest of this
section. 
\medskip

\noindent From \eqref{radon}, the phase function of $\mathcal{R}_1$ is 
$$\Phi_1(s,t,\vx_0,\alpha,\beta;\vx;\sigma)=\sigma\Psi_1(s,t,\vx_0,\alpha,\beta;\vx).$$
We now show that $\Phi_1$ is clean, non-degenerate and homogeneous in
$\sigma$ order 1, so that $\mathcal{A}$ satisfies the definition of
FIO (see definition \ref{FIOdef}). $\Phi_1$ is trivially homogeneous
order 1, since $\Psi_1$ does not depend on $\sigma$. Note also,
$\mathrm{d}_{s}\Phi_1=-\sigma\neq 0$, hence
$\mathrm{d}\Phi_1,\mathrm{d}_{\vy}\Phi_1\neq 0$. The apple surfaces
are smooth manifolds away from their singular points--the
 points which we do not consider. Hence $\Phi_1$ is clean. 

Let $\vxo =(x_0,y_0, z_0)$, then we will let $x'_{x_0}$ denote the
partial derivative of $x'=x'(\vx_0,\alpha,\beta;\vx)$ with respect to
$x_0$ and define the other partial derivatives of $x', y', z'$
analogously. Let
$$R^T=\begin{pmatrix}
\cos\alpha &\sin\alpha& 0 \\
-\sin\alpha\cos\beta &\cos\alpha\cos\beta& \sin\beta \\
\sin\alpha\sin\beta & -\cos\alpha\sin\beta & \cos\beta
\end{pmatrix}=\begin{pmatrix}r_1\\r_2\\r_3\end{pmatrix}$$
have rows $r_1,r_2,r_3$. Then, we have
\begin{equation}
\label{gradPhi}
\begin{split}
\nabla_{\vx_0}\Phi_1&=-2\sigma\left[\vx_T+\frac{t}{g}\paren{x'_{x_0}x'+y'_{x_0}y',x'_{y_0}x'+y'_{y_0}y',x'_{z_0}x'+y'_{z_0}y'}^T\right]\\
&=-2\sigma\paren{\vx_T+\frac{t}{g}\paren{
\nabla_{\vx_0}x', \nabla_{\vx_0}y'}\begin{pmatrix} r_1\\ r_2\end{pmatrix}\vx_T}\\
&=-2\sigma\paren{I-\frac{t}{g}A}\vx_T\\
&=-\sigma\nabla_{\vx}\Psi_1\\
&=-\mathrm{d}_{\vx}\Phi_1,
\end{split}
\end{equation}
where \bel{def:A}A=\paren{r_1^T,r_2^T}\begin{pmatrix} r_1\\
r_2\end{pmatrix}\ee is
symmetric, idempotent (i.e., $A^T=A$ and $A^2=A$) and , and $I$ is the
$3\times 3$ identity matrix. We have
\begin{equation}
\begin{split}
\text{det}\paren{I-\frac{t}{g}A}&=\text{det}\paren{I-\frac{t}{g}\paren{r_1^T,r_2^T}\begin{pmatrix} r_1\\ r_2\end{pmatrix}}\\
&=\text{det}\paren{I_{2\times 2}-\frac{t}{g}\begin{pmatrix} r_1\\ r_2\end{pmatrix}\paren{r_1^T,r_2^T}},\ \ (\text{by SDT})\\
&=\text{det}\paren{I_{2\times 2}-\frac{t}{g}I_{2\times 2}}\\
&=\paren{1-\frac{t}{g}}^2,
\end{split}
\end{equation}
which is zero if and only if $t=g$. Recall that we exclude the case $g=0$ since the singular points of
apples parameterized by $Y$ or $\tY$ do not meet $\oB$.
 
If $t\neq g$, then $\paren{I-\frac{t}{g}A}$ is invertible and
$\paren{I-\frac{t}{g}A}\vx_T=0\implies \vx_T=0$ but this would
mean that the center of the apple, $\vxo$, is on the apple
(equivalently, $t^2-s=0$). However, $s>t^2>0$ so this is not possible.

Now, we consider the case when $t=g$. Let $C_t$ denote the cylinder of
radius $t$ with axis of revolution $\{\vx_0+pR\veth : p\in\mathbb{R}\}$. If
$t=g$, then $\vx$ is in $ C_t$ and in
$\text{Null}(I-A)=\text{span}(r_1^T,r_2^T)$. Thus, if $t=g$, and
$(I-A)\vx_T=0$, then $\vx\in\{\vx_0+\text{span}(r_1^T,r_2^T)\}\cap
C_t$. If $\vx$ is in the critical set of $\Phi_1$ also (i.e., $\vx$
lies on the torus parameterized by $\vy$) then $s$ must be zero (i.e.,
the apple radius is zero), which we do not consider since
$\sqrt{s}>t>0$. Therefore, $\mathrm{d}_{\vx}\Phi_1\neq 0$ and $\Phi_1$
is nondegenerate. 

The amplitude of $\mathcal{A}$ is
\begin{equation}
\begin{split}
a_1(s,t,\vx_0,\alpha,\beta;\vx)&=\norm{\nabla_{\vx}\Psi_1(s,t,\vx_0,\alpha,\beta;\vx)}\\
&=2\paren{\vx_T^T\paren{I-\frac{t}{g}A}^T\paren{I-\frac{t}{g}A}\vx_T}^{\frac12}\\
\end{split}
\end{equation} by \eqref{gradPhi}.
By the arguments of the last paragraph we can show that $a_1$ is never zero, and hence $a_1$ is an elliptic symbol. $a$ is order zero since it is smooth, and does not depend on $\sigma$. Hence, $\mathcal{A}$ is an elliptic FIO order $O(\mathcal{A})=0+\frac{1}{2}-\frac{7}{2}=-2$.
\end{proof}

We now have our first main theorem which shows that $\mathcal{A}$ satisfies the semiglobal Bolker condition.

\begin{theorem}
\label{main1}
The left projection $\Pi^{(1)}_L$ of $\mathcal{A}$ is an injective immersion, and hence $\mathcal{A}$ satisfies the semiglobal Bolker condition from domain $\Ec'(B)$ to $\Dc'(Y)$.
\end{theorem}

%As noted in Remark \ref{rem:coords1} and the start of the proof
%of Theorem \ref{AFIO-1} we need to prove the theorem by parameterizing
%apples using the manifold without boundary $\tY$ from \eqref{tY map)},  the set $Y'$ from \eqref{def:Y}, and the map \eqref{tY map)}.

%To get local coordinates on $\tY$, we choose a vertical direction, $\ve_3$, and,
%from that vertical direction we establishes local coordinates, $(\alpha,\beta)$, on
%the set of direction axes, $\rpt$. Then, we define $Y'$ (see
%\eqref{def:Y'}) and the coordinate map \eqref{tY map)} using $(\alpha,\beta)$ for the direction axes. The remaining coordinates used in this study to parameterize spindle tori, namely $s$, $t$, and $\vx_0$,
%are global.

As the proof for Theorem \ref{main1} is long, we split the proof
into two subsections. We start with the immersion proof in the next
section, and present proof of injectivity in the following section. 

\subsubsection{$\Pi^{(1)}_L$ immersion proof}\label{sect:A immersion} Since being an
immersion is a local property , we can check this at an
arbitrary point $(\vy,\eta,\vx,\xi)$ in the canonical relation of
$\Ac$. Let $\ell$ be the direction axis of revolution of the apple parameterized by $\vy$.  Choose a unit vector $\ve_3$ such that $\ell$ is neither parallel nor perpendicular to $\ve_3$.  Choose unit vectors $\ve_1$ and $\ve_2$ so that $(\ve_1,\ve_2,\ve_3)$ makes up a right-hand coordinate system on $\rthree$.  Use this coordinate system on $\rthree$ to define the coordinate map \eqref{tY map)} and the set $Y'$ (see \eqref{def:Y'}).  Throughout this proof, the calculations are performed using
this coordinate system.

First we calculate $\Pi^{(1)}_L$. We have the derivatives
$$\mathrm{d}_t\Phi_1=2\sigma(t-g),\ \ \ \ \ \mathrm{d}_s\Phi_1=-\sigma,$$
\begin{equation}
\mathrm{d}_{\alpha}\Phi_1=\frac{2\sigma t}{g}\vx_T^T\paren{r_{1\alpha}^T,r_{2\alpha}^T}\begin{pmatrix} r_1\\ r_2\end{pmatrix}\vx_T,
\end{equation}
and
\begin{equation}
\mathrm{d}_{\beta}\Phi_1=\frac{2\sigma t}{g}\vx_T^T\paren{r_{1\beta}^T,r_{2\beta}^T}\begin{pmatrix} r_1\\ r_2\end{pmatrix}\vx_T,
\end{equation}
where $r_{i\alpha}$ is the component-wise partial derivative of $r_i$
with respect to $\alpha$ (similarly for $r_{i\beta}$). Let
\[H=\{\vx\in\mathbb{R}^3 \st
\Psi_1(t^2,t,\vx_0,\alpha,\beta;\vx)=0\}\] be the horn torus with
radius $t$ and axis of revolution $\vx_0+\rr R\ab\veth$, and let
$\tilde{H}=\{\vx\in\mathbb{R}^3 :
%\mathcal{T}
\Psi_1(t^2,t,\vx_0,\alpha,\beta;\vx)>0\}$ be the exterior of
$H$. Let $\mathcal{D}_1=\{(t,\vx_0,\alpha,\beta;\vx;\sigma) :
\vx\in\tilde{H}\}\times \mathbb{R}\backslash\{0\}$, then the map
\[\Dc_1\ni (t,\vx_0,\alpha,\beta;\vx;\sigma)\mapsto (s, t, \vx_0, \alpha, \beta
;\vx; \dd_s\Phi_1,\dd_t\Phi_1,
\dd_\beta\Phi_1;\nabla_{\vx_0}\Phi_1;\nabla_\vx\Phi_1)\in \Cc_1 \]
where $s = h-2tg$ gives local coordinates on the canonical relation for
$\Ac=\Rc_1$.

In these coordinates, the left projection $\Pi^{(1)}_L
:\mathcal{D}_1\to \Pi^{(1)}_L\paren{\mathcal{D}_1}$ of $\mathcal{A}$
is defined
\begin{equation}\label{def:Pi1}
\begin{split}
\Pi^{(1)}_L(\sigma;t,\alpha,\beta,\vxo;\vx)
=\Bigg(
&\overbrace{-\sigma}^{\mathrm{d}_s\Phi_1},
t,
\alpha,
\beta,
\vx_0,
\overbrace{-2\sigma\vx_T^T\paren{I-\frac{t}{g}A}}^{\nabla_{\vx_0}\Phi_1},
\overbrace{2\sigma(t-g)}^{\mathrm{d}_t\Phi_1},
\\
&\underbrace{\frac{2\sigma t}{g}\vx_T^T\paren{r_{1\alpha}^T,r_{2\alpha}^T}\begin{pmatrix} r_1\\ r_2\end{pmatrix}\vx_T}_{\mathrm{d}_{\alpha}\Phi_1},\underbrace{\frac{2\sigma t}{g}\vx_T^T\paren{r_{1\beta}^T,r_{2\beta}^T}\begin{pmatrix} r_1\\ r_2\end{pmatrix}\vx_T}_{\mathrm{d}_{\beta}\Phi_1},\underbrace{h-2tg}_{s}\Bigg),
\end{split}
\end{equation}
where we have highlighted the derivatives of $\Phi_1$ using under and
overbraces.  
Also, we have rearranged the variables in \eqref{def:Pi1} to correspond to the order used in calculating the Jacobian matrix of $\Pi^{(1)}_L$:
\begin{equation}\label{DPi1}
D\Pi^{(1)}_L=\kbordermatrix {&
\mathrm{d}\sigma,\mathrm{d}t,\mathrm{d}\alpha,\mathrm{d}\beta,
D_{\vxo} & \nabla_\vx \\
\dd_s\Phi_1,t,\alpha,\beta,{\vx_0} & I'_{7\times7} & \textbf{0}_{7\times3} \\
\nabla_{\vx_0}\Phi_1 & \cdot & D_\vx\paren{-2\sigma\vx_T^T\paren{I-\frac{t}{g}A}} \\
\mathrm{d}_t\Phi_1 & \cdot & \nabla_{\vx}\paren{2\sigma(t-g)}^T\\
\mathrm{d}_{\alpha}\Phi_1 & \cdot & \nabla_{\vx}\paren{\frac{2\sigma
t}{g}\vx_T^T\paren{r_{1\alpha}^T,r_{2\alpha}^T}\begin{pmatrix}
r_1\\ r_2\end{pmatrix}\vx_T}^{T}\\
\mathrm{d}_{\beta}\Phi_1 & \cdot & \nabla_{\vx}\paren{\frac{2\sigma
t}{g}\vx_T^T\paren{r_{1\beta}^T,r_{2\beta}^T}\begin{pmatrix} r_1\\
r_2\end{pmatrix}\vx_T}^{T}\\
s & \cdot & \nabla_{\vx}\paren{h-2tg}^{T}\\
},
\end{equation}
where $I'_{n\times n}$ is the $n\times n$ identity matrix but with the
first entry replaced by $-1$. Here we have highlighted the arguments
of $\Pi^{(1)}_L$ on the left-hand side of the matrix for
$D\Pi^{(1)}_L$, and the order of derivatives is indicated above
$D\Pi^{(1)}_L$. The terms corresponding to $\cdot$ in $D\Pi^{(1)}_L$
are not important for our calculations, as they will be multiplied by
zero in the calculation of the determinant of $D\Pi^{(1)}_L$. We now
find the derivatives in the right-hand column of
$D\Pi^{(1)}_L$ and show that $D\Pi^{(1)}_L$ is full rank.

Using the product rule and
$$\nabla_{\vx}\paren{\frac{1}{g}}=-\frac{1}{g^3}A\vx_T,$$
we can calculate the Jacobian matrix
\begin{equation}
\begin{split}
D_\vx\paren{-2\sigma\vx_T^T\paren{I-\frac{t}{g}A}}&=-2\sigma I+2\sigma\paren{-\frac{t}{g^3}A\vx_T\vx_T^TA^T+\frac{t}{g}A}\\
&=-2\sigma\paren{I-\frac{t}{g}A\paren{I-\frac{1}{g^2}\vx_T\vx_T^TA^T}}.
\end{split}
\end{equation}
Hence, using Sylvester's Determinant Theorem, it follows that
\begin{equation}
\label{J2}
\begin{split}
-\frac{1}{(2\sigma)^3}\text{det}\paren{D\paren{-2\sigma\vx_T^T\paren{I-\frac{t}{g}A}}}&=\text{det}\paren{I-\frac{t}{g}A\paren{I-\frac{1}{g^2}\vx_T\vx_T^TA^T}}\\
&=\text{det}\paren{I-\frac{t}{g}\paren{I-\frac{1}{g^2}\vx_T\vx_T^TA^T}A}\\
&=\text{det}\paren{I-\frac{t}{g}\paren{I-\frac{1}{g^2}\vx_T\vx_T^T}A}\\
&=\text{det}\paren{I+C}\\
&=1+\text{tr}(C)+\frac{1}{2}\paren{(\text{tr}(C))^2-\text{tr}(C^2)}+\text{det}(C),
\end{split}
\end{equation}
where we use SDT in the second step to reverse the matrix multiplication order, and the fact that $A$ is symmetric idempotent in the third step to get $A^TA=A^2=A$. Here $C=-\frac{t}{g}\paren{I-B}A$, where $B=\frac{1}{g^2}\vx_T\vx_T^T$.

We now simplify \eqref{J2}. First, we have the identities
\begin{equation}
\begin{split}
\text{tr}(A)&=\text{tr}\paren{\paren{r_1^T,r_2^T}\begin{pmatrix} r_1\\ r_2\end{pmatrix}}=\text{tr}\paren{\begin{pmatrix} r_1\\ r_2\end{pmatrix}\paren{r_1^T,r_2^T}}=\text{tr}\paren{I_{2\times 2}}=2,
\end{split}
\end{equation}
\begin{equation}
\begin{split}
\text{tr}(AB)&=\frac{1}{g^2}\text{tr}\paren{\vx_T\vx_T^TA}=\frac{1}{g^2}\text{tr}\paren{\vx_T^TA\vx_T}=\frac{\vx_T^TA\vx_T}{g^2}=1,
\end{split}
\end{equation}
noting that $\vx_T^TA\vx_T=g^2$,
$$\text{tr}\paren{BABA}=\frac{1}{g^4}\text{tr}\paren{\vx_T(\vx_T^TA\vx_T)\vx_T^TA}=\frac{\vx_T^TA\vx_T}{g^4}\text{tr}\paren{\vx_T\vx_T^TA}=\frac{(\vx_T^TA\vx_T)^2}{g^4}=1,$$
and
$$\text{tr}(ABA)=\text{tr}(BA^2)=\text{tr}(BA)=1.$$
Now
$$
\text{tr}(C)=-\frac{t}{g}\text{tr}(A-BA)=-\frac{t}{g}(\text{tr}(A)-\text{tr}(BA))=-\frac{t}{g},$$
and 
\[
\begin{aligned}
\frac{1}{2}\paren{\tr(C))^2-\tr(C^2)}&=\frac{1}{2}
\bparen{\tr(C)^2 -\frac{t^2}{g^2}\paren{\tr(BABA)-\tr(ABA) -\tr(BAA) +\tr(A^2)}}
\\& =\frac{t^2}{2g^2}\bparen{1-(1 -1-1+2) }=0
\end{aligned}
\]
and
$$\text{det}(C)=-\frac{t^3}{g^3}\text{det}(A)\text{det}(I-B)=0,$$
since $\text{det}(A)=0$. Indeed $0\neq r_3^T\in\text{Null}(A)$. Putting this together, we have
$$\text{det}\paren{D\paren{-2\sigma\vx_T^T\paren{I-\frac{t}{g}A}}}=-(2\sigma)^3\paren{1-\frac{t}{g}},$$
which is zero if and only if $t=g$. Hence, in the case when $t\neq
g$, $D\Pi^{(1)}_L$ has full rank and $\Pi^{(1)}_L$ is an immersion.

We now consider the case when $t=g$. In this case, $\vx=R\vx_C+\vx_0$,
where $\vx_C=(t\cos\theta,t\sin\theta,z)^T$, for some
$\theta\in[0,2\pi]$ and $z\in \rr$. That is, $\vx$ lies on the
cylinder of radius $t$, with axis of revolution $\{\vx_0+pR\veth : p\in
\mathbb{R}\}$.

% \tred{But $\vx$ also lies on the apple
% $T_1(s,t,\vxo,\alpha,\beta)$ and so $z'^2 = s$ where we recall that
% $(x',y',z') = R^T\vx_T$. }
% \tc{I had hoped to use the above red text, but I don't think it helps.
%  If you don't see a way to use it, we can delete it.}

Under the assumption $t=g$, we show that  the submatrix
\begin{equation}M=
\kbordermatrix {& \nabla_\vx \\
\mathrm{d}_t\Phi_1  & \nabla_{\vx}\paren{2\sigma(t-g)}^T\\
\mathrm{d}_{\alpha}\Phi_1  & \nabla_{\vx}\paren{\frac{2\sigma t}{g}\vx_T^T\paren{r_{1\alpha}^T,r_{2\alpha}^T}\begin{pmatrix} r_1\\ r_2\end{pmatrix}\vx_T}^T\\
s & \nabla_{\vx}\paren{h-2tg}^T\\
}
\end{equation}
of $D\Pi^{(1)}_L$ is invertible. Using the product rule and
$$\nabla_{\vx}\paren{g^i}=\frac{i}{2}g^{i-2}\nabla_{\vx}\paren{x'^2+y'^2}=ig^{i-2}A\vx_T,$$
for $i\in\mathbb{Z}$, we have
\begin{equation}\label{M expression}M=
\kbordermatrix {& \nabla_\vx \\
\mathrm{d}_t\Phi_1  & -\frac{2\sigma}{g}(A\vx_T)^T\\
\mathrm{d}_{\alpha}\Phi_1  & 2\sigma\left[\frac{2t}{g}\paren{r_{1\alpha}^T,r_{2\alpha}^T}\begin{pmatrix} r_1\\ r_2\end{pmatrix}\vx_T-\frac{t}{g^3}\paren{\vx_T^T\paren{r_{1\alpha}^T,r_{2\alpha}^T}\begin{pmatrix} r_1\\ r_2\end{pmatrix}\vx_T}A\vx_T\right]^T\\
s & 2\paren{\vx_T-\frac{t}{g}A\vx_T}^T
}.
\end{equation}
Substituting $\vx=R\vx_C+\vx_0$ and $t=g$, we have
\begin{equation}
\label{equ2.20}
M=
\begin{pmatrix}
-\frac{2\sigma}{t}\left[(r_1^T,r_2^T,0)\vx_C\right]^T\\
2\sigma\left[2\paren{r_{1\alpha}^T,r_{2\alpha}^T,0}\vx_C+\frac{1}{t^2}(tz\sin\beta\cos\theta)(r_1^T,r_2^T,0)\vx_C\right]^T\\
2\paren{R\vx_C-(r_1^T,r_2^T,0)\vx_C}^T,
\end{pmatrix}
\end{equation}
where
\begin{equation}
\label{3.21}
%\hspace{-1cm}
\begin{split}
\vx_T^T\paren{r_{1\alpha}^T,r_{2\alpha}^T}\begin{pmatrix} r_1\\ r_2\end{pmatrix}\vx_T&=\vx_C^TR^T\paren{r_{1\alpha}^T,r_{2\alpha}^T}\begin{pmatrix} r_1\\ r_2\end{pmatrix}R\vx_C\\
&=\vx_C^TR^T\paren{r_{1\alpha}^T,r_{2\alpha}^T,0}\vx_C\\
%&=\vx_C^T\begin{pmatrix}
%\cos\alpha &-\sin\alpha\cos\beta& \sin\alpha\sin\beta \\
%\sin\alpha &\cos\alpha\cos\beta& -\cos\alpha\sin\beta \\
%0 &\sin\beta& \cos\beta
%\end{pmatrix}^T\begin{pmatrix}
%-\sin\alpha &-\cos\alpha\cos\beta& 0 \\
%\cos\alpha &\sin\alpha\cos\beta& 0 \\
%0 & 0 & 0
%\end{pmatrix}\vx_C\\
&=\vx_C^T\begin{pmatrix}
0 & 0 & 0 \\
\cos\beta & 0 & 0 \\
-\sin\beta & 0 & 0
\end{pmatrix}\vx_C\\
&=- tz\sin\beta\cos\theta.
\end{split}
\end{equation}
Here \eqref{3.21} shows the calculations for the new scalar term in brackets on the second row of $M$ in \eqref{equ2.20}. We have
\begin{equation}MR=
\begin{pmatrix}
-2\sigma(\cos\theta,\sin\theta,0)\\
2\sigma\left[2\paren{-t\cos\beta\sin\theta,t\cos\beta\cos\theta,-t\sin\beta\cos\theta}+z\sin\beta\cos\theta(\cos\theta,\sin\theta,0)\right]\\
2\paren{0,0,z}
\end{pmatrix}.
\end{equation}
Thus
\begin{equation}
\begin{split}
\text{det}(M)&=\text{det}(MR)\\
&=8\sigma^2\begin{pmatrix}0\\0\\ z\end{pmatrix}
\cdot\left[-\begin{pmatrix}\cos\theta\\\sin\theta\\ 0\end{pmatrix}\times\left[2\begin{pmatrix}-t\cos\beta\sin\theta\\t\cos\beta\cos\theta\\ -t\sin\beta\cos\theta\end{pmatrix}+z\sin\beta\cos\theta\begin{pmatrix}\cos\theta\\\sin\theta\\ 0\end{pmatrix}\right]\right]\\
&=-16\sigma^2zt\begin{pmatrix}0\\0\\ 1\end{pmatrix}
\cdot\left[\begin{pmatrix}\cos\theta\\\sin\theta\\ 0\end{pmatrix}\times\begin{pmatrix}-\cos\beta\sin\theta\\\cos\beta\cos\theta\\ -\sin\beta\cos\theta\end{pmatrix}\right]\\
&=-16\sigma^2zt\cos\beta.
\end{split}
\end{equation}
Recall that $\beta\in(0,\pi/2)$ by definition of $Y'$, and $t>0$. The
case $z=0$ corresponds to $s=0$, i.e., degenerate tori which have
radius zero and collapse into a circle of radius $t$ passing through the
center of the apple tube. We do not consider  degenerate tori.
Hence $M$, and thus $D\Pi^{(1)}_L$, have full rank and $\Pi^{(1)}_L$
is an immersion.

\subsubsection{$\Pi^{(1)}_L$ injectivity proof}  Injectivity of
$\Pi_L^{(1)}$ is a local property in the target space. To determine if
$\Pi_L^{(1)}$ is injective, we take an arbitrary point $(\vy,\eta)\in T^*(\tY)$
and see if it has more than one preimage. Specifically, we choose
$\vy\in \tY$ and take local coordinates \eqref{tY map)} so that $\vy$
is in the image of $Y'$ (i.e., the axis of the spindle torus parametrized by $\vy$ is neither vertical nor horizontal). Then, we analyze $\Pi_L^{(1)}$ using these
coordinates.

Let $\vx_1, \vx_2\in B$ be such that 
$$\Pi^{(1)}_L(t,\vx_0,\alpha,\beta;\vx_1;\sigma)=\Pi^{(1)}_L(t,\vx_0,\alpha,\beta;\vx_2;\sigma).$$
Then
$$2\sigma(t-g_1)=2\sigma(t-g_2)\implies g_1=g_2=g,$$
where $g_1,g_2$ correspond to the inputs $\vx_1,\vx_2$. We now consider two cases, namely $t=g$ and $t\neq g$.\\
\\
\textbf{Case 1: $t\neq g$}\\
We have
$$\paren{I-\frac{t}{g}A}\paren{\vx_1-\vx_0}=\paren{I-\frac{t}{g}A}\paren{\vx_2-\vx_0}.$$
Thus, $\Pi^{(1)}_L$ is injective if $\paren{I-\frac{t}{g}A}$ is
invertible. Following similar arguments to those used in Theorem
\ref{AFIO-1}, we have
\begin{equation}
\label{3.30}
\begin{split}
\text{det}\paren{I-\frac{t}{g}A}&=\text{det}\paren{I-\frac{t}{g}\paren{r_1^T,r_2^T}\begin{pmatrix} r_1\\ r_2\end{pmatrix}}\\
&=\text{det}\paren{I_{2\times 2}-\frac{t}{g}\begin{pmatrix} r_1\\ r_2\end{pmatrix}\paren{r_1^T,r_2^T}}\\
&=\text{det}\paren{I_{2\times 2}-\frac{t}{g}I_{2\times 2}}\\
&=\paren{1-\frac{t}{g}}^2\neq 0.
\end{split}
\end{equation}
Therefore, $\vx_1=\vx_2$ and $\Pi^{(1)}_L$ is injective. Note we have used SDT in the second step of \eqref{3.30} to reverse the matrix multiplication order inside the determinant.\\
\\
\textbf{Case 2: $t=g$}\\
In this case, $\vx_j=R\vx^{(j)}_C+\vx_0$, for $j=1,2$, where
$\vx^{(j)}_C=(t\cos\theta_j,t\sin\theta_j,z_j)^T$ and
$\theta_j\in[0,2\pi]$. That is, the $\vx_j$ lie on the cylinder,
radius $t$, with axis of revolution $R\veth+\vx_0$.

Using \eqref{gradPhi} and \eqref{def:A}, we have
\begin{equation}
\begin{split}
\nabla_{\vx_0}\Phi_1(t,\vx_0,\alpha,\beta;\vx_j;\sigma)&=-2\sigma\left[R\vx^{(j)}_C-AR\vx^{(j)}_C\right]\\
&=-2\sigma\left[(r_1^T,r_2^T,r_3^T)\vx^{(j)}_C-(r_1^T,r_2^T,0)\vx^{(j)}_C\right]\\
&=-2\sigma  z_j r_3^T.
\end{split}
\end{equation}
Hence $z_1=z_2=z\neq 0$, since $r_3^T\neq \textbf{0}$, and we do not
consider the case $z=\sqrt{s}=0$ (i.e., a degenerate torus). 
%\jc{\tred{We may need a smooth cutoff for the circle through the center of the torus tube, even though we never integrate there. What I mean is $z=0$ is in the domain of $\Pi_L$ currently, but $z=0$ is also bounded away from the surfaces of integration. This seems similar to our last paper where $\Pi_L$ had problems on vertical lines $x=x_0$, and we needed a cutoff, even though the curves of integration never intersected such lines.}}
Note, the $\vx_j$ are constrained also to lie on the apple parameterized by $(s,t,\vx_0,\alpha,\beta)$, which, in the $t=g$ case, implies $z_j=\sqrt{s}$. Now,
\begin{equation}
%\hspace{-1cm}
\begin{split}
\mathrm{d}_{\alpha}\Phi_1(t,\vx_0,\alpha,\beta;\vx_j;\sigma)&=2\sigma (\vx^{(j)}_C)^TR^T\paren{r_{1\alpha}^T,r_{2\alpha}^T,0}\vx^{(j)}_C\\
%&=2\sigma (\vx^{(j)}_C)^T\begin{pmatrix}
%\cos\alpha &-\sin\alpha\cos\beta& \sin\alpha\sin\beta \\
%\sin\alpha &\cos\alpha\cos\beta& -\cos\alpha\sin\beta \\
%0 &\sin\beta& \cos\beta
%\end{pmatrix}^T\begin{pmatrix}
%-\sin\alpha &-\cos\alpha\cos\beta& 0 \\
%\cos\alpha &\sin\alpha\cos\beta& 0 \\
%0 & 0 & 0
%\end{pmatrix}\vx^{(j)}_C\\
&=2\sigma (\vx^{(j)}_C)^T\begin{pmatrix}
0 & -\cos\beta & 0 \\
\cos\beta & 0 & 0 \\
-\sin\beta & 0 & 0
\end{pmatrix}\vx^{(j)}_C\\
&=-2\sigma tz\cos\theta_j\sin\beta,
\end{split}
\end{equation}
and
\begin{equation}
\begin{split}
\mathrm{d}_{\beta}\Phi_1(t,\vx_0,\alpha,\beta;\vx_j;\sigma)&=2\sigma (\vx^{(j)}_C)^TR^T\paren{r_{1\beta}^T,r_{2\beta}^T,0}\vx^{(j)}_C\\
%&=2\sigma (\vx^{(j)}_C)^T\begin{pmatrix}
%\cos\alpha &-\sin\alpha\cos\beta& \sin\alpha\sin\beta \\
%\sin\alpha &\cos\alpha\cos\beta& -\cos\alpha\sin\beta \\
%0 &\sin\beta& \cos\beta
%\end{pmatrix}^T\begin{pmatrix}
%0 & \sin\alpha\sin\beta& 0 \\
%0 & -\cos\alpha\sin\beta & 0 \\
%0 & \cos\beta & 0
%\end{pmatrix}\vx^{(j)}_C\\
&=2\sigma (\vx^{(j)}_C)^T\begin{pmatrix}
0 & 0 & 0 \\
0 & 0 & 0 \\
0 & 1 & 0
\end{pmatrix}\vx^{(j)}_C\\
&=-2\sigma tz\sin\theta_j.
\end{split}
\end{equation}
Therefore, $(\cos\theta_1,\sin\theta_1)=(\cos\theta_2,\sin\theta_2)$ and
$\Pi^{(j)}_L$ is injective. Recall that $\beta\in(0,\pi/2)$, because
$(s,t,\vxo,\alpha,\beta)\in Y'$, and so $
\sin\beta>0$. 
\\

\noindent This completes the proof of Theorem \ref{main1}.

% \tc{Above, I would like to see if there is a row of $\nabla_\vx
% \paren{\nabla_{\vxo}\Phi_1}$ that can show $\theta_1 = \theta_2$ when
% $\beta = 0$. I feel like it should happen since the problem is
% rotation invariant and the problem when $\beta =0$ and our map
% $(\alpha,\beta)\mapsto R^T(\alpha,\beta)$ does not
% provide coordinates. Since this problem is rotation invariant, one
% should be able to use different coordinates that are good near
% $\veth$. This is not a problem for the immersion proof since that
% doesn't depend on $\alpha$ or $\beta$, and I'll add some text to this
% effect. For the injectivity proof, I hope we can use one line of
% $D_\vx\paren{\nabla_\vxo \Phi_1}$ for the case when $t=g$.}

\subsection{Microlocal properties of $\mathcal{L}$; the $j=2$ case}
Here we discuss the microlocal properties of $\mathcal{L}$ in a
similar way to the $j=1$ case. First, we prove that $\mathcal{L}$ is an elliptic FIO order $-2$.

\begin{theorem}
\label{AFIO-2}
The lemon transform $\mathcal{L}=\mathcal{R}_2$ of \eqref{radon} is an elliptic FIO order $-2$ from domain $\Ec'(B)$ to $\Dc'(Y)$.
\end{theorem}

\begin{proof}  As in Theorem \ref{AFIO-1}, we choose local coordinates \eqref{tY map)} on $\tY$ and use these local coordinates in our calculations.

From \eqref{radon}, the phase function of $\mathcal{R}_2$ is
$$\Phi_2(s,t,\vx_0,\alpha,\beta;\vx;\sigma)=\sigma\Psi_2(s,t,\vx_0,\alpha,\beta;\vx).$$
We now show that $\Phi_2$ is clean, non-degenerate and homogeneous in
$\sigma$ order 1, to show that $\mathcal{L}$ satisfies the definition
of FIO (see definition \ref{FIOdef}). $\Phi_2$ is trivially
homogeneous order 1. $\mathrm{d}_{s}\Phi_2=-\sigma\neq 0$, and hence
$\mathrm{d}\Phi_2,\mathrm{d}_{\vy}\Phi_2\neq 0$. The lemon surfaces
are smooth manifolds away from their singular points, which we do not
consider by the definition of $Y$, \eqref{def:Y}. Hence $\Phi_2$ is clean. 

Using similar calculations to those of \eqref{gradPhi}, we have
\begin{equation}
\begin{split}
\nabla_{\vx_0}\Phi_2&=-2\sigma\paren{I+\frac{t}{g}A}\vx_T\\
&=-\sigma\nabla_{\vx}\Psi_2\\
&=-\mathrm{d}_{\vx}\Phi_2.
\end{split}
\end{equation}
Also
\begin{equation}
\begin{split}
\text{det}\paren{I+\frac{t}{g}A}&=\text{det}\paren{I+\frac{t}{g}\paren{r_1^T,r_2^T}\begin{pmatrix} r_1\\ r_2\end{pmatrix}}\\
&=\text{det}\paren{I_{2\times 2}+\frac{t}{g}\begin{pmatrix} r_1\\ r_2\end{pmatrix}\paren{r_1^T,r_2^T}},\ \text{by (SDT)}\\
&=\text{det}\paren{I_{2\times 2}+\frac{t}{g}I_{2\times 2}}\\
&=\paren{1+\frac{t}{g}}^2>0.
\end{split}
\end{equation}
Thus $\mathrm{d}_{\vx}\Phi_2$ is zero if and only if $\vx_T=0$, which we do not consider. Hence $\Phi_2$ is nondegenerate.

The amplitude is
\begin{equation}
\begin{split}
a_2(s,t,\vx_0,\alpha,\beta;\vx)&=\norm{\nabla_{\vx}\Psi_2(s,t,\vx_0,\alpha,\beta;\vx)}\\
&=2\paren{\vx_T^T\paren{I+\frac{t}{g}A}^T\paren{I+\frac{t}{g}A}\vx_T}^{\frac{1}{2}}.
\end{split}
\end{equation}
$a_2$ is smooth, and independent of $\sigma$, and hence $a_2$ is a
symbol order zero. $a_2>0$ since $I+\frac{t}{g}A$ is invertible, and
hence $\paren{I+\frac{t}{g}A}^T\paren{I+\frac{t}{g}A}$ is positive
definite. Therefore $\mathcal{L}$ is an elliptic FIO order
$O(\mathcal{L})=0+\frac{1}{2}-\frac{7}{2}=-2$. 

Recall that the spindle
tori in $\tY$ do not have singular points in $\oB$. Therefore, $g$ is
never zero and the symbol is defined for functions are supported in
$B$.\end{proof}

\noindent We now have our second main theorem which shows that $\mathcal{L}$ satisfies the semiglobal Bolker condition.

\begin{theorem}
\label{main2}
The left projection $\Pi^{(2)}_L$ of $\mathcal{L}$ is an injective immersion, and hence $\mathcal{L}$ satisfies the semiglobal Bolker condition.
\end{theorem}

We now proceed in a similar fashion to the proof of Theorem \ref{main2}, i.e., we split the proof into two subsections. We start with the immersion proof in the next section, and prove injectivity in the following section.

\subsubsection{$\Pi^{(2)}_L$ immersion proof}
We choose a point $(\vy,\eta,\vx,\xi)$ in the canonical relation of $\Lc$ and  choose coordinates as for the apple transform in section \ref{sect:A immersion} so the spindle torus axis is neither vertical nor horizontal.

The left projection of $\mathcal{R}_2$ is
\begin{equation}
\label{Pi2}
\begin{split}
\Pi^{(2)}_L(t,\vx_0,\alpha,\beta;\vx;\sigma)=&\Bigg(\overbrace{h+2tg}^{s},t,\vx_0,\alpha,\beta,\overbrace{-\sigma}^{\mathrm{d}_s\Phi_2},\overbrace{-2\sigma\vx_T^T\paren{I+\frac{t}{g}A}}^{\nabla_{\vx_0}\Phi_2},\\
&\,\,\,\,\underbrace{2\sigma(t+g)}_{\mathrm{d}_t\Phi_2},\underbrace{-\frac{2\sigma t}{g}\vx_T^T\paren{r_{1\alpha}^T,r_{2\alpha}^T}\begin{pmatrix} r_1\\ r_2\end{pmatrix}\vx_T}_{\mathrm{d}_{\alpha}\Phi_2},\underbrace{-\frac{2\sigma t}{g}\vx_T^T\paren{r_{1\beta}^T,r_{2\beta}^T}\begin{pmatrix} r_1\\ r_2\end{pmatrix}\vx_T}_{\mathrm{d}_{\beta}\Phi_2}\Bigg).
\end{split}
\end{equation}
The proof is analogous to the $j=1$ case, and a little easier,
so we will go over the main points.  

We calculate $D \Pi_L^{(2)}$ and  just consider the
rows corresponding to $D_{\vx} \paren{\nabla_{\vxo} \Phi_2}$.  We can
show, in a similar way to the $j=1$ case,
$$\text{det}\paren{D\paren{-2\sigma\vx_T^T\paren{I+\frac{t}{g}A}}}=-(2\sigma)^3\paren{1+\frac{t}{g}},$$
which is never zero. Therefore, these rows of $D \Pi_L^{(2)}$
have full rank $3$. Hence $D\Pi^{(2)}_L$ has full rank and
$\Pi^{(2)}_L$ is an immersion. 

\subsubsection{$\Pi^{(2)}_L$ injectivity proof}  
To prove injectivity we proceed similarly to Theorem \ref{AFIO-1}, i.e., we take an arbitrary point $(\vy,\eta)\in T^*(\tY)$
and determine whether it has more than one preimage under $\Pi_L^{(2)}$. We choose coordinates on $\rthree$ so that the axis of the lemon parameterized by $\vy$ is neither vertical nor horizontal.

Let $\vx_1, \vx_2\in B$ be such that 
$$\Pi^{(2)}_L(t,\vx_0,\alpha,\beta;\vx_1;\sigma)=\Pi^{(2)}_L(t,\vx_0,\alpha,\beta;\vx_2;\sigma)=(\vy,\eta).$$
Then
$$2\sigma(t+g_1)=2\sigma(t+g_2)\implies g_1=g_2=g,$$
where $g_1,g_2$ correspond to the inputs $\vx_1,\vx_2$. 

Focusing on the $\nabla_\vxo \Phi_2$ terms in the image of $\Pi_L^{(2)}$ (see \eqref{Pi2}), we have
$$\paren{I+\frac{t}{g}A}\paren{\vx_1-\vx_0}=\paren{I+\frac{t}{g}A}\paren{\vx_2-\vx_0}.$$
Thus, $\Pi^{(2)}_L$ is injective if $\paren{I+\frac{t}{g}A}$ is invertible. We have
\begin{equation}
\begin{split}
\text{det}\paren{I+\frac{t}{g}A}&=\paren{1+\frac{t}{g}}^2>0.
\end{split}
\end{equation}
Hence, $\Pi^{(2)}_L$ is injective.\\
\\
This concludes the proof of Theorem \ref{main2}.

\begin{corollary}
\label{corr1} Let $\alpha\in[0,2\pi]$, and $\beta\in[0,\pi/2]$ be fixed.
Let $f\in \Ec'(B)$ and $(s,t,\vxo)$ chosen so the singular points of the spindle torus parameterized by $(s,t,\vxo,\alpha,\beta)$ are disjoint from $\overline{B}$.
Then the Radon transform 
$$\mathcal{L}_Tf(s,t,\vx_0)=\mathcal{L}f(s,t,\vx_0,\alpha,\beta),$$
which defines the integrals of $f$ over a 5-D set of translated
lemons, satisfies the semi\-glo\-bal Bolker condition. 

The analogous restriction for the apple transform
$$\mathcal{A}_Tf(s,t,\vx_0)=\mathcal{A}f(s,t,\vx_0,\alpha,\beta),$$
however, does not satisfy the semiglobal Bolker condition.
\end{corollary}

\begin{proof}
This follows immediately from Theorems \ref{main1} and \ref{main2}.
The left projection of $\mathcal{A}_T$ has Jacobian which drops rank
on the cylinder $t=g$, and thus there are artifacts in the
reconstruction which occur along rings which are the intersections of
apples and cylinders, radius $t$, with the same axis of revolution.

Regarding $\mathcal{L}$, we require only a 3-D translation of the
lemons, and the radial variables ($s$ and $t$), in order for the
semiglobal Bolker condition to be satisfied. The rotations induced by $\alpha,\beta$ are not
needed in the proof of Theorem \ref{main2}. 
\end{proof}

\subsubsection{Discussion}
In \cite{webberholman}, lemon transforms are analyzed, but only
rotations and changes in radius of the lemons are considered (i.e.,
$\vx_0=0$ is fixed, and $\alpha$, $\beta$, and $s,t$ vary). The
authors prove that the left projection drops rank, and show that there
are artifacts in (unfiltered) backprojection, and Landweber image
reconstructions. With knowledge of seven-dimensional lemon integral data, however, we would not expect to see artifacts due to rank deficiencies in the reconstruction. In fact, five-dimensional lemon integral data is sufficient to show the Bolker condition is satisfied, as is shown by Corollary \ref{corr1}. 

With regards to $\mathcal{A}$, the full
seven-dimensional data is needed in the proof of Theorem \ref{main1}
to show that the Bolker condition is satisfied. As noted in Corollary \ref{corr1},
there is issue with the translated apples on their intersections with
cylinders, radius $t$, which share the same axis of revolution. Such
issues can be addressed by including the 2-D rotation induced by
$\alpha$ and $\beta$. Thus, with knowledge of seven-dimensional apple integral data, we would not expect to see artifacts due to microlocal properties in the
reconstruction.\\
\\
\noindent In the following section, we consider 3-D subsets of apple and lemon surfaces which have practical motivations in CST.
In the case of the apple transform, we discover artifacts which occur
at apple-cylinder intersections, and are thus consistent with the
results of Theorem \ref{main1}.

\section{Practical geometry in CST}
\label{CST}

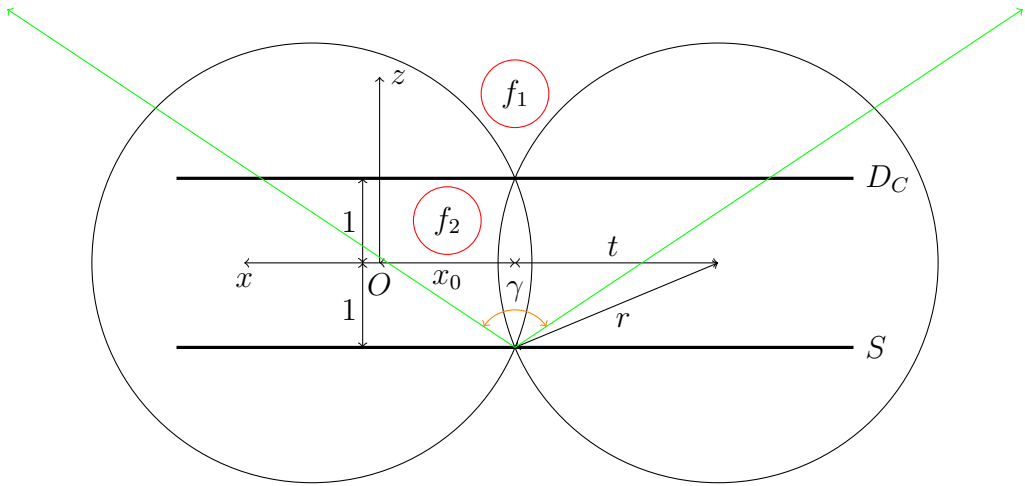
\begin{figure}[!h]
\centering
\begin{tikzpicture}[scale=4.5]
%\draw [dashed] (-1,0.15)--(0.76,0.15)node[right] {$D_A$ at $\{z=2-r_m\}$};%origianlly y=0
\draw [very thick] (-1,1)--(1,1)node[right] {$S$};
\draw [<->] (-0.1-0.35,1)--(-0.1-0.35,1.25);
\node at (-0.14-0.35,1.1125) {$1$};
\draw [<->] (-0.1-0.35,1.25)--(-0.1-0.35,1.5);
\node at (-0.14-0.35,1.3725) {$1$};
%\draw [<->] (-0.1-0.35,0.75)--(-0.1-0.35,1);
%\node at (-0.14-0.35,0.8725) {$1$};
\draw [very thick] (-1,1.5)--(1,1.5)node[right] {$D_C$};
%\draw [<->] (-0.4,1.55)--(1,1.55);
%\node at (0.3,1.58) {$a$};
\draw [->] (-0.4,1.25)--(-0.4,1.8)node[right] {$z$};
\draw [->] (-0.4,1.25)node[below] {$O$}--(-1.2+0.4,1.25)node[below] {$x$};
\draw [<->] (-0.4,1.25)--(0,1.25);
\node at (-0.2,1.2) {$x_0$};
\draw (0.6,1.25) circle (0.65);
\draw (-0.6,1.25) circle (0.65);
\draw [<->] (0.5+0.1,1.25)--(0,1);
\node at (0.24+0.08,1.08) {$r$};
\draw [<->] (0.5+0.1,1.25)--(0,1.25);
\node at (0.24+0.05,1.3) {$t$};
%\node at (0.04,1.75) {$x_2$};
%\node at (1.04,0.45+0.25) {$x_1$};
%\node at (-0.03,0.47+0.25) {$O$};
%\node at (0.32,0.62) {$\vx$};
%\draw [->]  (0.37,0.64)--(0.32,0.52)node[left] {$\vxi$};
%\draw [->]  (0.37,0.64)--(0.42,0.52);
\coordinate (origo1) at (0.37,0.64);
\coordinate (pivot1) at (0.42,0.52);
\coordinate (bob1) at (0.32,0.52);
%\draw pic[fill=orange, <->,"$C_T$", angle eccentricity=1.9] {angle = bob1--origo1--pivot1};
%\draw  (0.37,0.64)--(0.6,1.25);
%\node at (-0.6,2) {$C_1(s,x_0)$};
%\node at (0.6,2) {$C_2(s,x_0)$};
%\draw  [dashed] (0.37,0.64)--(0.37,1.25);
\coordinate (origo) at (0.37,0.64);
\coordinate (pivot) at (0.37,1.25);
\coordinate (bob) at (0.6,1.25);
%\draw pic[draw=orange, <->,"$\beta$", angle eccentricity=2] {angle = bob--origo--pivot};
%\draw [->] (-0.6,1.5)--(0.6,0);
%\draw [->] (0.6,1.5)--(-0.6,0);
%\coordinate (origo1) at (0,0.75);
%\coordinate (pivot1) at (0.6,0);
%\coordinate (bob1) at (0.6,0.75);
%\draw pic[draw=orange, <->,"$\alpha_m$", angle eccentricity=1.7] {angle = pivot1--origo1--bob1};
\draw [green,->] (0,1)--(1.5,2);
\draw [green,->] (0,1)--(-1.5,2);
\coordinate (origo) at (1.5,2);
\coordinate (pivot) at (0,1);
\coordinate (bob) at (-1.5,2);
\draw pic[draw=orange, <->,"$\gamma$", angle eccentricity=1.5] {angle = origo--pivot--bob};
\draw [red] (-0.2,1.375) circle (0.1);
\node at (-0.2,1.375) {$f_2$};
\draw [red] (0,1.75) circle (0.1);
\node at (0,1.75) {$f_1$};
\end{tikzpicture}
\caption{Two-dimensional cross-section of parallel CST geometry. $S$ and $D_C$ denote the source and detector rows, which  are parallel and one unit away from the origin. The $y$ axis is perpendicular to the page.  Here we show a 2-D cross section of a spindle torus with center $\vx_0$, and $t>0$ is the distance from $\vx_0$ to the center of the torus tube, as before. The torus has radius $r=\sqrt{t^2+1}>0$.
}
\label{fig1.4}
\end{figure}

In this section we consider the machine geometry of \cite{webber2019compton,WebberQuinto2020II}, which has practical applications in airport baggage screening. We present a microlocal analysis of the apple transform, first introduced in \cite{webber2019compton}, and its lemon transform analog. Specifically, we consider the machine geometry of figure \ref{fig1.4}.

%=\sqrt{s}

 The diagram illustrates an X-ray scanner comprised of a line segment of sources ($S$), which emit X-rays in the direction of a parallel line segment of detectors ($D_C$). The photons are then Compton scattered and measured by the detectors on $D_C$. Meanwhile, the target, $f$, is translated out of the page (i.e., in the $y$ axis direction) on a conveyor belt. We consider two possibilities for the location of $f$ here, namely within the half space $\{z>1\}$ (i.e., above $D_C$), and within the band $\{-1<z<1\}$ (i.e. between $S$ and $D_C$). Examples of these two possible locations for $f$ are illustrated by $f_1$ and $f_2$ in figure \ref{fig1.4}, respectively, where $f_1$ is integrated over apples and $f_2$ over lemons. The source is cone beam with opening angle $\gamma$. We set $\gamma=\pi$ so that photons are everywhere on $\{z>-1\}$.   
 See \cite{webber2019compton,WebberQuinto2020II}, for more details on the applications to airport baggage screening and CST, more generally. In total, the data is three-dimensional, and is comprised of a 2-D translation and a 1-D radial variable. In this section, $R=I$ is fixed (i.e., there is no rotation of the apples or lemons), $z_0=0$ (i.e., the translation is in the $(x,y)$ plane), and $s$ and $t$ satisfy the relation $s=t^2+1$ (see figure \ref{fig1.4}).
With this in mind we define the restricted apple and  lemon transforms 
 %$\mathcal{A} : L_c^2(\{z>1\})\to \mathcal{L}\paren{L_c^2(\{z>0\})}$
\bel{def:AoLo}\begin{gathered}\mathcal{A}_0f(p,x_0,y_0)=\mathcal{A}f\paren{\frac{p}{4}+1,\sqrt{\frac{p}{4}},(x_0,y_0,0)^T,0,0},\\
\text{and}\\
%$\mathcal{L} : L_c^2(B)\to \mathcal{L}\paren{L_c^2(B)}$
\mathcal{L}_0f(p,x_0,y_0)=\mathcal{L}f\paren{\frac{p}{4}+1,\sqrt{\frac{p}{4}},(x_0,y_0,0)^T,0,0},\\
\text{where}\ \ \ p=4t^2.\end{gathered}\ee The variable $p$ is introduced in this section to
simplify the  calculations. 
%In this case $f\in L_c^2(B)$ when integrating over lemons, as before. For the apple integrals, $f\in L_c^2(\{z>1\})$ is supported above the detector array, as is the case in (\tred{reference}).

\begin{proposition}
The restricted apple and  lemon transforms can be written
\begin{equation}\label{def:Ao Lo}
\mathcal{A}_0f=\mathcal{T}_1f\ \text{for}\ f\in L_c^2(\{z>1\}),
\ \ \ \ 
\mathcal{L}_0f=\mathcal{T}_2f,\ \text{for}\ f\in L_c^2(\{-1<z<1\})
\end{equation}
where
\begin{equation}
\label{Tradon}
\begin{split}
\mathcal{T}_jf(p,x_0,y_0)&=\int_{X_j}\norm{\nabla_{\vx}\Psi(p,x_0,y_0;\vx)}\delta\paren{\Psi(p,x_0,y_0;\vx)}f(\vx)\mathrm{d}\vx\\
&=\int_{-\infty}^{\infty}\int_{X_j}\norm{\nabla_{\vx}\Psi(p,x_0,y_0;\vx)}e^{\sigma
\Psi(p,x_0,y_0;\vx)}f(\vx)\mathrm{d}\vx\mathrm{d}\sigma,
\end{split}
\end{equation}
for $j=1,2$, where we now define 
\begin{gather}X_1=\sparen{\vx\in \rthree\st z>1}, \quad X_2=\sparen{\vx\in\mathbb{R}^3 :
(x-x_0)^2+(y-y_0)^2+z^2<1},\notag\\
\Psi(p,x_0,y_0;\vx)=p-\frac{h^2}{g}\quad \text{where}\label{def:Psi}\\
g=g(x_0,y_0; \vx)=(x-x_0)^2+(y-y_0)^2,\quad h=h(x_0,y_0;
\vx)=\norm{\vx_T}^2-1,\notag\\
\vx_T=(x-x_0,y-y_0,z)\notag\end{gather}

\end{proposition}
\noindent Note that the functions $\Psi, g, h, \vx_T$ are adapted from  section \ref{radon section} for our geometry.

\begin{proof}
A torus centered at the origin with axis of rotation $z$ is described
implicitly by the equation
\begin{equation}
\label{torus}
\paren{\norm{\vx}^2+t^2-s}^2=4t^2(x^2+y^2).
\end{equation}
Hence, the defining equation for the tori of interest, which are translated by $(x_0,y_0)$ in the $(x,y)$ plane (as depicted in figure \ref{fig1.4}), and satisfy $t^2-s=-1$, becomes
\begin{equation}
\Psi(p,x_0,y_0;\vx)=p-\frac{\paren{(x-x_0)^2+(y-y_0)^2+z^2-1}^2}{(x-x_0)^2+(y-y_0)^2}=
p-\frac{h^2(x_0,y_0;\vx)}{g(x_0,y_0;\vx)}
\end{equation}
%where $p=4t^2$.
 Thus, when the integration is restricted to $X_1$, $\mathcal{T}_1f$
 defines the integrals of $f$ over a 3-D set of translated apples
 whose singular points lie on $\{z=1\}$ and $\{z=-1\}$. $\mathcal{T}_2f$ defines integrals of $f$ over lemons in the same way when the integration is restricted to $X_2$. 
 
Note that the functions in $L^2_c\paren{\sparen{z>1}}$ and
$L^2_c\paren{\sparen{-1<z<1}}$ the domains of $\Ac_0$ and $\Lc_0$,
respectively, in \eqref{def:Ao Lo} are zero near the singular points
of the spindle tori (which satisfy $z=\pm 1$). Hence the surface
measure on spindle tori for  $\mathcal{T}_j$ is defined on the support of
$f$.\end{proof}

\noindent We now show that the $\mathcal{T}_j$ are elliptic FIO order $-1$.

%Note that, due to the support restrictions on functions in $f$, $g>0$, i.e., we
%stay away from the singular points of the torus at $z=\pm 1$.

\begin{theorem}
The Radon transforms $\mathcal{T}_j$, for $j=1,2$, are elliptic FIO order $-1$ from domain $\Ec'(\sparen{z>1})$ for $j=1$ and from domain $\Ec'(\sparen{-1<z<1}$ for $j=2$.
\end{theorem}

\begin{proof}
For the proofs in this section, it will be convenient to define the
function 
\begin{equation}
u(x_0,y_0;\vx)=\frac{h(x_0,y_0; \vx)}{g(x_0,y_0;\vx)}.
\end{equation}

The phase function of $\mathcal{T}_j$, for $j=1,2$, is 
\bel{def:Phi}\Phi(p,x_0,y_0;\vx;\sigma)=\sigma\Psi(p,x_0,y_0;\vx)=\sigma\paren{p-hu},\ee
by \eqref{Tradon},
and the amplitude is
\begin{equation}\label{def:a}
a(p,x_0,y_0; \vx)=\norm{\nabla_{\vx}\Psi(p,x_0,y_0;
\vx)}=2u\sqrt{g(u-2)^2+4z^2}.
\end{equation}

%Both are undefined on the rotation axis of the spindle torus, when
%$g=0$, however this is only a technicality.  

The phase \eqref{def:Phi} and the amplitude \eqref{def:a} are
undefined when $g=0$, that is on the rotation axis of each spindle
torus--when $(x,y) =\xoyo$. To get around this, we use a smooth cutoff
near the spindle torus axis as in \cite[Lemma
3.3]{webber2021microlocal} to smoothly set the symbol to zero near the
spindle torus axis. Note that the points at which the cutoff is not
smooth, the singular points of the spindle torus (on $z=\pm 1$), are
not in either domain in \eqref{def:Ao Lo}. This cutoff makes the
amplitude defined and smooth everywhere. Note that the phase is smooth
on a neighborhood of the canonical relations of our transforms and the
cutoff on the symbol can be used to make it smooth everywhere. Similarly, for $\Lc_0$, the integral in \eqref{Tradon}  for each $\xoyo$ is over the open disk, $X_2$ so as to integrate only over the lemon, not the part of the apple in $\sparen{-1<z<1}$. One constructs a smooth function of $(p,\xo,\yo;\vx)$  that is equal to $1$ in a neighborhood in $(0,\infty)\times\rtwo\times \sparen{-1<z<1}$ of the lemon parameterized by $(p,\xo,\yo)$ and equal to $0$ in a neighborhood of the corresponding apple.

The phase \eqref{def:Phi} is trivially clean and homogeneous in
$\sigma$ order 1. In addition, $\mathrm{d}_p\Phi=\sigma\neq 0$, and
\[
\norm{\mathrm{d}_{\vx}\Phi}=2\sigma u\sqrt{g(u-2)^2+4z^2}.
\]
$u=0$ does not occur since $p>0$. When $j=1$, the domain of integration is such that $z>1$, and hence $\mathcal{T}_2$ has nondegenerate phase. When $j=2$, $-1<z<1$, and 
$$u-2=\frac{1}{g}\left[z^2-(x-x_0)^2-(y-y_0)^2-1\right]<0.$$
Hence $\mathcal{T}_1$ has nondegenerate phase. 

By the same arguments as for the phase, the amplitude \eqref{def:a} is
positive on the apple and on the lemon, the manifolds of integration
of $\Ac_0$ and $\Lc_0$. By using the cutoff, $a$ is smooth.
Furthermore, $a$ does not depend on $\sigma$, so $a$ is an elliptic
symbol order zero, and thus $\mathcal{T}_j$, for $j=1,2$, is an
elliptic FIO order $O(\mathcal{T}_j)=0+\frac{1}{2}-\frac{3}{2}=-1$.
\end{proof}

We now have our third main theorem which provides conditions such that the $\mathcal{T}_j$ satisfy the semiglobal Bolker condition.

\begin{theorem}
\label{main3}Global coordinates on the canonical relation of
$\Tc_j$ are given by $(x_0,y_0;\vx;\sigma)$, as in \eqref{tY map)}
but with $\ab = (0,0)$ and $p =
h^2(x_0,y_0;\vx)/g(x_0,y_0;\vx)$.   Let
\bel{def:D1}D_1=\mathbb{R}^2\times\{\vx\in\mathbb{R}^3 :
z>1\}\times \mathbb{R}\backslash\{0\}\ee and
\bel{def:D2}D_2=\{(x_0,y_0;\vx)\in\mathbb{R}^5 :
z\in (0,1)\}\times
\mathbb{R}\backslash\{0\}.\ee  
These sets define global coordinates on the appropriate canonical relation.  

The left projection $\Pi^{(1)}_L : D_1\to
\Pi^{(1)}_L\paren{D_1}$ of $\mathcal{T}_1$ is an injective immersion
under the constraint that \bel{u cond} u-2 >0, \ \
\text{equivalently}\ \ z^2-1 >(x-\xo)^2 - (y-\yo)^2.\ee 

The left projection $\Pi^{(2)}_L : D_2\to \Pi^{(2)}_L\paren{D_2}$ of
$\mathcal{T}_2$ is an injective immersion.
\end{theorem}

%\tc{I think you can change $D_2=\{(x_0,y_0; \vx)\in\mathbb{R}^5 :
%z\in (0,1)\}\times
%\mathbb{R}\backslash\{0\}$}

\begin{remark}\label{rem:D1}The requirement in Theorem \ref{main3}
that functions are supported in the half-space $z>0$ (or equivalently $z<0$) is natural
because the apples and lemons are symmetric about $z=0$. Therefore,
our transforms integrate odd functions in $z$ to zero and
singularities for $z<0$ can cancel singularities for $z>0$.

We point out that \eqref{u cond} puts restrictions on the
support of functions and the sets of $\xoyo$ for which $\Ac_0$
satisfies the Bolker condition. Define the set \bel{def:H} H\xoyo =
\sparen{(x,y,z)\st z>1,\ z^2-1>(x-\xo)^2-(y-\yo)^2}.\ee Let $f\in
L^2_c(\sparen{z>1})$. If $\supp(f)$ is so large that it is not contained in $H\xoyo$ for any $\xoyo$, then one cannot apply
Theorem \ref{main3} to $f$.

Now, assume $\xoyo\in \rtwo$ and $K$ is a compact subset of
$z>1$ such that $K\subset H\xoyo$. Then by compactness of $K$, there
is an open neighborhood, $U$, of $\xoyo$ such that $K\subset
H(x_1,y_1)$ for all $(x_1,y_1)\in U$. Theorem \ref{main3} can be
applied to the local problem for $\Ac_0$ for functions supported in
$K$ and centers $\xoyo$ in $U$.
\end{remark}

\begin{proof}
The points $(\xo,\yo;\vx;\sigma)$ are  coordinates on the
canonical relation of $\Tc_j$ for the same reason as  \eqref{tY map)}
give local coordinates on $\tY$. However, here they are global  
coordinates because $\ab=(0,0)$ is fixed.

Let $u_1=u-1$ and $u_2=u-2$. Then, the left projection of $\mathcal{T}_j$, for $j=1,2$, is
\begin{equation}
\Pi^{(j)}_L(\sigma;x_0,y_0;\vx)=\paren{\overbrace{\sigma}^{\mathrm{d}_p\Phi},x_0,y_0,\overbrace{uh}^{p},\overbrace{-2\sigma uu_2(x-x_0)}^{\mathrm{d}_{x_0}\Phi}, \overbrace{-2\sigma uu_2(y-y_0)}^{\mathrm{d}_{y_0}\Phi}},
\end{equation}
where we have rearranged the variables to correspond to the order used in calculating the Jacobian
\bel{Jacobian 1}D\Pi^{(j)}_L=\begin{pmatrix} I_{3\times 3} & 0_{3\times 3}\\\cdot & M_2\end{pmatrix},\ee
where, using
$$\mathrm{d}_x(uu_2)=2u_xu_1=-\frac{4}{g}(x-x_0)u_1^2,$$
and
$$\mathrm{d}_y(uu_2)=2u_yu_1=-\frac{4}{g}(y-y_0)u_1^2,$$
we have
\begin{equation}\label{M2}
\begin{split}
M_2=\kbordermatrix {& \mathrm{d}x & \mathrm{d}y & \mathrm{d}z \\
p &{ -2u_2(x-x_0)u} & {-2u_2(y-y_0)u} &{4zu}\\
\mathrm{d}_{x_0}\Phi &{- \frac{2\sigma}{g} (-4(x-x_0)^2u_1^2+hu_2)} & {\frac{8\sigma}{g} (x-x_0)(y-y_0)u_1^2} & {\frac{-8\sigma}{g}(x-x_0)zu_1}\\
\mathrm{d}_{y_0}\Phi &{ \frac{8\sigma}{g} (x-x_0)(y-y_0)u_1^2} & {\frac{-2\sigma}{g}(-4(y-y_0)^2u_1^2+hu_2)}
& {-\frac{8\sigma}{g}(y-y_0)z u_1}}.
\end{split}
\end{equation}
The determinant of $M_2$ is hence
\begin{equation}
\text{det}(M_2)=\frac{16z\sigma^2}{g^2}u\times\text{det}(M_3),
\end{equation}
where
\begin{equation}
{M_3=
\begin{pmatrix} -(x-x_0)u_2 & -(y-y_0)u_2 & 1\\
 4(x-x_0)^2u_1^2-hu_2 & 4(x-x_0)(y-y_0)u_1^2 &-2 (x-x_0)u_1\\
 4(x-x_0)(y-y_0)u_1^2 & 4(y-y_0)^2u_1^2-hu_2 & -2(y-y_0)u_1\end{pmatrix}}.
\end{equation}

%\noindent Letting $M_3=\begin{pmatrix} r_1\\ r_2\\ r_3\end{pmatrix}$ %have rows $r_1,r_2,r_3$, we have
%\begin{equation}
%\begin{split}
%r_2\times r_3&=\\
%&\paren{-hu_1u_2(x-x_0),-hu_1u_2(y-y_0),-4ghu_1^2u_2+h^2u_2^2}.
%\end{split}
%\end{equation}
%Therefore
%\begin{equation}
%\begin{split}
%\text{det}(M_3)&=r_1\cdot(r_2\times r_3)\\
%&=2hu_1u_2^2\left[(x-x_0)^2+(y-y_0)^2\right]-4ghu_1^2+h^2u_2^2\\
%&=2ghu_1u_2^2-4ghu_1^2+h^2u_2^2\\
%&=-2ghu_1u_2\paren{2u_1-u_2}+h^2u_2^2\\
%&=-2ughu_1u_2+h^2u_2^2\\
%&=-h^2u_2\paren{-u_2+2u_1},\ \ \ (\text{using $ug=h$})\\
%&=-h^2uu_2.\\
%\end{split}
%\end{equation}

\noindent A straightforward calculation shows that $\text{det}(M_3)=-h^2uu_2$, and hence
\begin{equation}
\text{det}D\Pi^{(j)}_L=-16z\sigma^2u^4(u-2).
\end{equation}
Therefore, $\text{det}D\Pi^{(j)}_L=0$ if and only if $z=0$, $u=0$ or $u=2$.
On $D_1$, $u,z>0$ and hence the left projection of $\mathcal{T}_1$
drops rank if and only if $u=2$. By condition \eqref{u cond},  $u>2$, and hence
$\text{det}D\Pi^{(1)}_L$ is nonzero and $\Pi_L^{(1)}$ is  an immersion. 
On $D_2$, $u_2,u<0$ and
hence the left projection of $\mathcal{T}_{2}$ drops rank if and only if
$z=0$, which we do not consider by assumption that $z>0$. Hence
$\text{det}D\Pi^{(2)}_L$ is an immersion.

Now onto injectivity. First, we consider the case $j=1$. Let $\vx_1=(x_1,y_1,z_1)$ and $\vx_2=(x_2,y_2,z_2)$ be such that 
$$\Pi^{(1)}_L(x_0,y_0; \vx_1;\sigma)=\Pi^{(1)}_L(x_0,y_0;\vx_2;\sigma),$$
and let $v=u(x_0,y_0; \vx_1)$, and $w= u(x_0,y_0; \vx_2)$.
Let  $v_1=v-1$, $v_2=v-2$,   $w_1=w-1$, and $w_2=w-2$.
%Further let $h_j=h(\vx_j,x_0,y_0)$, and $g_j=g(\vx_j,x_0,y_0)$, for $j=1,2$. Then
\begin{equation}\label{equal}
\begin{split}
\paren{vh_1,(x_1-x_0)vv_2,(y_1-y_0)vv_2}=\paren{wh_2,(x_2-x_0)ww_2,(y_2-y_0)ww_2},
\end{split}
\end{equation}
where
$h_j=h(x_0,y_0; \vx_j)$, and $g_j=g(x_0,y_0; \vx_j)$, for $j=1,2$. It follows that
\begin{equation}
\label{inj}
\begin{split}
&\left[(x_1-x_0)^2+(y_1-y_0)^2\right]v^2v_2^2=\left[(x_2-x_0)^2+(y_2-y_0)^2\right]w^2w_2^2\\
&\implies g_1v^2v_2^2=g_2w^2w_2^2\\
&\implies vh_1v_2^2=wh_2w_2^2\\
&\implies v_2^2=w_2^2,\ \ \ (\text{note $vh_1=wh_2=p>0$}).
\end{split}
\end{equation}
%Let $f_1-2>0$ and $f_2-2<0$. Then $(x_1-x_0)=-(x_2-x_0)$ and $(y_1-y_0)=-(y_2-y_0)$ (note $f_j\neq 0,2$, so we can divide by $f_j(f_j-2)$), which implies that $f_1-2<0$ and $f_2-2>0$ and we have a contradiction. To explain this further, $\vx_1$ and $\vx_2$ cannot be on either side of the hyperboloid $\{f=2\}$, and satisfy $\vx_1=R\vx_2$, where $R$ is a rotation about $z$, since $\{f=2\}$ is a surface of rotation about $z$. Hence $f_1-2$ and $f_2-2$ must be the same sign. 
Under our assumption that $v_2,w_2>0$, it follows that $v_2=w_2$, so $v=w$ and  $h_1=h_2$. Using \eqref{equal} again, we see that $v=w\implies (x_1,y_1)=(x_2,y_2)$ (note $v,w\neq
0,2$), and so $z_1^2= z_2^2$ since $h_1=h_2$.
%Here, $z^2$ is determined from the definition of $h$. 
On $D_1$, $z>1$, so $z_1=z_2$ and $\Pi^{(1)}_L$ is
thus an injective immersion.

%\tc{Note:  it seems that one could also %show $\Pi_1^{(1)}$ satisfies Bolker for %$z>1$ and $u_2<0$, i.e.,  between $z=1$ %and the hyperboloid $H\xoyo$.}

On $D_2$, $v_2,w_2<0$, and hence by \eqref{inj} and the
previous arguments $\Pi^{(2)}_L(x_0,y_0;
\vx_1;\sigma)=\Pi^{(2)}_L(x_0,y_0;\vx_2;\sigma)$ implies that $(x_1,y_1)=(x_2,y_2)$ and $z_1^2= z_2^2$. Also, $z\in(0,1)$, which implies $z_1=z_2$, and thus $\Pi^{(2)}_L$ is an injective immersion.
\end{proof}

\subsection{Discussion of the artifacts}
In this section, we discuss the restrictions imposed on the function support and left projection domain in Theorem \ref{main3} needed to show that the semiglobal Bolker condition is satisfied, and address the artifacts that occur when such constrains are lifted. 
\begin{figure}[!h]
\centering
\begin{tikzpicture}[scale=1.4]
\draw [very thick] (-4,-1)--(4,-1)node[right] {$D_C$};
\draw [very thick] (-4,1.5-0.5)--(4,1.5-0.5)node[right] {$S$};
\draw [->] (-1,0)--(-1,0.75)node[right] {$z$};
\draw [->] (-1,0)node[below] {$O$}--(-1.75,0)node[below] {$x$};
\draw [<->] (-1,0)--(0,0);
\node at (-0.5,-0.15) {$x_0$};
\draw (1.73,0) circle (2);
\draw (-1.73,0) circle (2);
\draw [<->] (-2,1)--(-2,-1);
\node at (-2.1,0) {2};
\draw [<->] (1.73,0)--(1.73,2);
\node at (1.83,1.25) {$r$};
\draw [blue,domain=-2.1:2.1] plot(\x, {sqrt(1+pow(\x,2))});
\draw [blue,domain=-2.1:2.1] plot(\x, {-sqrt(1+pow(\x,2))});
\draw [->] (2.4,2.4)node[right] {$\{u-2=0\}$}--(2.1,2.32);
\end{tikzpicture}
\caption{Intersection of a torus and $\{u-2=0\}$ shown as a 2-D cross-section. A 2-D cross-section of $\{u-2=0\}$ is drawn in blue and intersects the translated torus at $z=\pm r$.}
\label{fig2}
\end{figure}
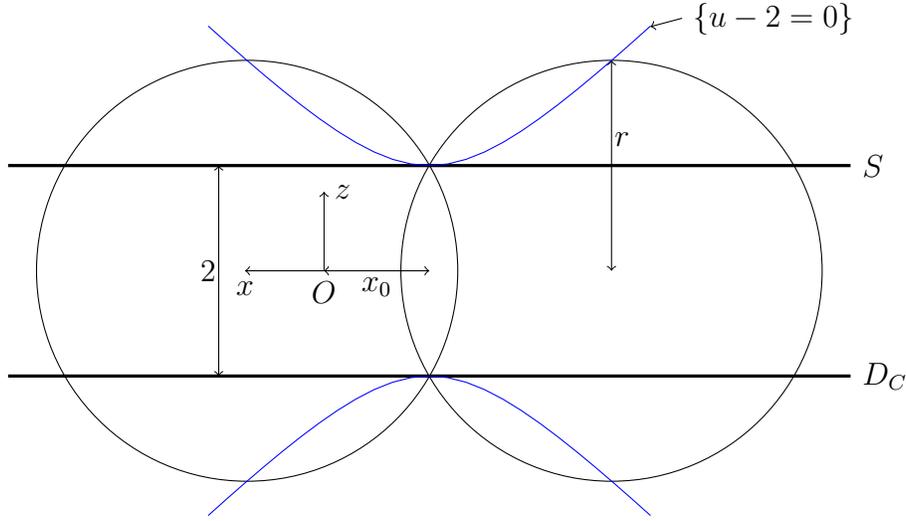

%\tc{FYI:  Figure 5 shows that the bad points on the spindle torus are the top and bottom points (also the points at $z=0$).  This is because you can't wiggle $\xoyo$ and $p$ enough to get nearby cotangent directions at those points.  Namely, you can't wiggle perpendicular to the plane of figure 5.}

In the proof of Theorem \ref{main3}, the constraint $z\in (0,1)$ was needed to show that $\Pi^{(2)}_L$ is an injective immersion. Without loss of generality, we could replace the $z\in (0,1)$ constraint with $z\in (-1,0)$ and the proof would follow in the same way, by symmetry of the apples and lemons about $z=0$. If the function support is not restricted in this way, and $z$ takes values on the full range $-1<z<1$, then $\Pi^{(2)}_L$ is noninjective. For example, $\Pi^{(2)}_L((x,y,z);\sigma,x_0,y_0)=\Pi^{(2)}_L((x,y,-z);\sigma,x_0,y_0)$, for any $-1<z<1$, and thus there are artifacts which consist of reflections in the $(x,y)$ plane. When $z=0$, specifically, $D\Pi^{(2)}_L$ drops rank and $\Pi^{(2)}_L$ is nonimmersive.  $\Pi^{(1)}_L$ also suffers the same noninjectivity concerns if $z$ takes both signs. However, it does not practically make sense for the function to be supported on both sides of the $(x,y)$ plane when integrating over apples. Indeed, the cone-beam direction in figure \ref{fig1.4} is such that there are no photons on $\{z<-1\}$, and thus we assume $f$ is supported on $\{z>1\}$ when integrating over apples (as is done also in \cite{WebberQuinto2020II}). 
%The Radon transforms $\mathcal{T}_j$, for $j=1,2$, include in their null space the antisymmetric functions in $z$ (e.g., $f=z$ is in $\text{Null}(\mathcal{T}_j)$ for $j=1,2$). Hence, without such positivity (or negativity) constraints on the $z$ integration bounds, there is a nontrivial null space. 

Regarding $\Pi^{(1)}_L$, $u>2$ was assumed a-priori in Theorem \ref{main3}, and discussed also in Remark \ref{rem:D1}, in order to show that $\Pi^{(1)}_L$ satisfies the Bolker condition. Without such restrictions, in particular when $u=2$, $D\Pi^{(1)}_L$ drops rank and there artifacts which occur along rings at the top and bottom of the apple surface. Specifically, when $u=2$, $\vx$ lies on the two-sided hyperboloid, described implicitly by
\begin{equation}
\label{hyp}
z^2-(x-x_0)^2-(y-y_0)^2-1=0.
\end{equation}

The intersection of the apple and the surface defined by \eqref{hyp} occurs when $z=\pm \sqrt{t^2+1}= r$, i.e., along the rings at the top and bottom of the apple. See figure \ref{fig2}, where we have shown a 2-D cross section of the intersecting apple and hyperboloid surfaces. In Corollary \ref{corr1}, we showed that $\mathcal{A}_T$ did not satisfy the Bolker condition. Specifically, the left projection of $\mathcal{A}_T$  drops rank for $\vx$ on $\{t=g\}$, namely the cylinders radius $t$, with axis of revolution $\{\vx_0+\nu R\veth : \nu\in\mathbb{R}\}$, i.e., the axis of revolution of the apple surface. The apple and $\{t=g\}$ intersect on rings at the top and bottom of the apple which are the same intersection points as those shown in figure \ref{fig2}. Thus, our results are consistent with the findings of Corollary \ref{corr1}. Specifically, when the degrees of freedom in our data includes translation, e.g., the full 3-D translation of $\mathcal{A}_Tf$, or the 2-D translation of $\mathcal{A}_0f$, there is a consistency in the artifact locations. 
%\tred{In practice, for $f$ with singularities in the $z$ direction (i.e., in directions normal to the rings at the top and bottom of the apple), we can expect to see a blurring effect over such rings in a reconstruction from translated apple integral data (e.g., $\mathcal{A}_Tf$ or $\mathcal{A}_0f$).}

%\tc{Are you using the same intuition I used in the FYI I gave at the start of the section about how, at those points, one cannot perturb the normal in all directions?  Even if that was your motivation, I think more needs to be said or it can be for another article in which you do the reconstructions. }
%Such artifacts are typically highlighted in Filtered Back-Projection (FBP) and Landweber type reconstructions, as is, for example, seen  in the reconstructions of \cite{webberholman}.

%\jc{\tred{Given where the blue curves and tori intersect (i.e., at the top and bottom of the torus where the normal directions are parallel to the $z$ axis), I think we're going to see artifacts on planes parallel to the $xy$ plane whenever the function has a singularity parallel to the $z$ axis. Could we come up with a way to filter such artifacts in an FBP?}}
%So $\Pi_L$ is an immersion away from $\{f=2\}$ and $\{z=0\}$. Going back to equation \eqref{amplitude}, we can see that $a$ is never zero, since both $f-2$ and $z$ cannot be zero at the same time (note that also $g,f\neq 0$). Hence, $\mathcal{R}$ is an elliptic FIO. 

\subsection{How to remove artifacts with machine design}
In this section, we discuss possible modifications to the machine design of figure \ref{fig1.4}, so that the conditions of Theorem \ref{main3} are met, and thus we do not have to contend with the types of artifacts discussed in the previous section. 

\begin{figure}[!h]
\centering
\begin{tikzpicture}[scale=1.4]
\draw [very thick] (-4,-1)--(4,-1)node[right] {$D_C\  \{z=-1\}$};
\draw [very thick] (-4,1.5-0.5)--(4,1.5-0.5)node[right] {$S\ \{z=1\}$};
\draw [red,dashed] (-4,1.5-0.5+0.3)--(4,1.5-0.5+0.3)node[right] {$\{z=1+\epsilon\}$};
\draw [red,dashed] (-4,0)--(4,0)node[right] {$\{z=0\}$};
\draw [->] (-1,0)--(-1,0.75)node[right] {$z$};
\draw [->] (-1,0)node[below] {$O$}--(-1.75,0)node[below] {$x$};
\draw [<->] (-1,0)--(0,0);
\node at (-0.5,-0.15) {$x_0$};
\draw (1.73,0) circle (2);
\draw (-1.73,0) circle (2);
%\draw [<->] (1.73,0)--(1.73,2);
%\node at (1.83,1.2) {$r$};
\draw [blue,domain=-2.1:2.1] plot(\x, {sqrt(1+pow(\x,2))});
%\draw [blue,domain=-2.1:2.1] plot(\x, {-sqrt(1+pow(\x,2))});
\draw [->] (2.4,2.4)node[right] {$\{z=\sqrt{1+x^2}\}$}--(2.1,2.32);
\draw [red] (1.3,0.5) circle (0.3);
\node at (1.3,0.5) {$f_2$};
\draw [red] (0,2) circle (0.3);
\node at (0,2) {$f_1$};
\draw [green,->] (0,-1)--(1.3,2.5);
\draw [green,->] (0,-1)--(-1.3,2.5);
\coordinate (origo) at (1.3,2.5);
\coordinate (pivot) at (0,-1);
\coordinate (bob) at (-1.3,2.5);
\draw pic[draw=orange, <->,"$\gamma$", angle eccentricity=1.5] {angle = origo--pivot--bob};
\end{tikzpicture}
\caption{Restrictions to source cone beam angle (backscatter) and object support (forward scatter) so that Bolker is satisfied.}
\label{fig3}
\end{figure}
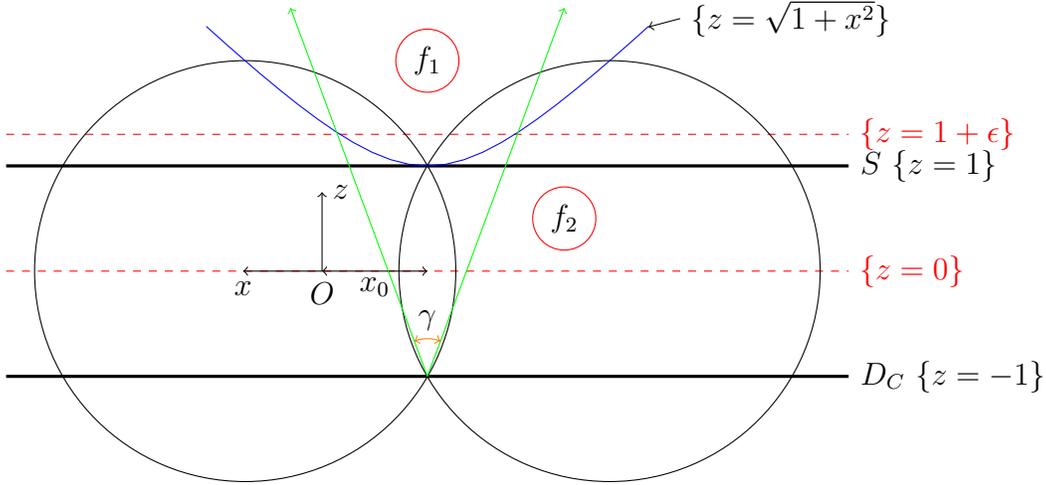

When using forward scattered photons for imaging, whose intensity is modeled by the lemon transform, we need only restrict the support of $f$ to $\{0<z<1\}$ or $\{-1<z<0\}$. See figure \ref{fig3} for an example $f$ with such support, in particular the location of $f_2$. In this case, the conditions of Theorem \ref{main3} are satisfied and the lemon transform satisfies the Bolker condition. Practically speaking, such support restrictions can be achieved by re-positioning the scanning target (e.g., the airport luggage) to be strictly above or below the $(x,y)$ plane. For example, we could construct the conveyor belt to lie on $\{z=0\}$ (highlighted by a red dashed line in figure \ref{fig3}) and place the scanning target (with height less than 1) on top of the conveyor, to ensure the conditions of Theorem \ref{main3} are met.

%\tc{In this  discussion below one needs to move $\xoyo$ in an open neighborhood, not just a line to get 3 dimension of data and for Bolker to hold, so the conveyor belt above needs to shimmy. \newline Also, what about not collimating the source but just  moving $\xoyo$ in an open set $U\subset \rtwo$ so $\supp(f)$ stays inside the hyperboloid $H\xoyo$ for all $\xoyo\in U$?   Equivalently, one would fix $\xoyo$ and just move the object around on the conveyor so it stays inside this fixed $H\xoyo$.  This would be equivalent to what I just described, it would be easier to do, and you don't collimate data.}

Regarding $\mathcal{A}_0$, i.e., when backscattered photons are used for imaging, the object is compactly supported on $\{z>1\}$. To ensure the conditions of Theorem \ref{main3} are met, we propose to further restrict the support of $f$ to $\{z>1+\epsilon\}$, for some $\epsilon>0$. See $f_1$ in figure \ref{fig3} for an example $f$ with such support. In practice, this would mean placing the conveyor belt on $\{z=1+\epsilon\}$ (shown as a red dashed line in figure \ref{fig3}), with $f$ on top of the conveyor. With such restrictions on the support of $f$, we can choose the cone-beam angle $\gamma$ (as shown in figure \ref{fig3}) so that no scatter occurs on the surface $\{u-2=0\}$, and $u-2>0$. Note that we have removed the bottom half of $\{u-2=0\}$ in figure \ref{fig3}, since $f_1$ is supported on $\{z>1+\epsilon\}$. To restrict the scatter exclusively to $u-2>0$, we can write $\gamma$ explicitly as
\begin{equation}
\gamma=\min\paren{2\tan^{-1}\frac{\sqrt{(1+\epsilon)^2-1}}{2+\epsilon},90^{\circ}}.
\end{equation}
Note, $\gamma$ must be less than or equal to $90^{\circ}$ since if $\gamma>90^{\circ}$, the line through $(x_0,-1)$ and $(x_0\pm\sqrt{(1+\epsilon)^2-1},1+\epsilon)$ (the right$\backslash$left hand green line of figure \ref{fig3}) has gradient less than 1, and hence intersects the blue curve of figure \ref{fig3} for large enough $x$ (i.e., scatter could occur on or below $\{u-2=0\}$). Note that the blue curve in figure \ref{fig3} (i.e., $\{z=\sqrt{1+x^2}\}$) has max gradient 1, so the green lines of figure \ref{fig3} must have gradient greater than or equal to 1 to ensure that $u-2>0$ and the conditions of Theorem \ref{main3} are satisfied. In practice, such restrictions on $\gamma$ would mean there is less signal, due to the smaller cone-beam and less photons, and hence the data would become more noisy. So, while we can address the microlocal artifacts by restricting $\gamma$, this would in turn increase the noise level. Thus, there is a trade of to consider here, i.e., do we want higher Signal to Noise Ratio (SNR) and more artifacts, or less SNR with less artifacts? We leave such practical concerns for future work. 

\section{Conclusions and further work} In this paper, we presented a
novel microlocal analysis of seven-dimensional apple and lemon Radon
transforms, which have applications in CST. The goal of this work was
to consider a best case scenario in CST, in terms of data
dimensionality. The literature \cite{webberholman, me2, rigaud20183d,
webber2019compton, rigaud20213d, rigaud2021reconstruction,
cebeiro2021three} considers exclusively Radon transforms which define
the integrals of a function over three-dimensional sets of apple or
lemon surfaces. In these works, artifacts are present in the
reconstruction due to data limitations, and regularization strategies
are used to combat the artifacts. Here, we considered a case when a
full seven-dimensional set of apple and lemon integrals are known. Our
main theorems, namely Theorems \ref{main1} and \ref{main2}), prove that
the apple and lemon transforms are elliptic FIO, order 2, which
satisfy the Bolker condition. 

%\tred{Following our main theorems, we
%proved that the normal operators of the apple and  lemon transforms are
%elliptic PDO, and thus have a bounded inverse in Sobolev norms order
%4. Therefore, with sufficient data in CST, image artifacts due to microlocal properties are absent
%and there exist modest bounds on the level of noise amplification. More
%specifically, in terms of noise amplification, an inverse apple or 
%lemon transform operation is equivalent to taking two derivatives (i.e., a
%mildly ill-posed inversion).}

In addition, we investigated the microlocal properties of apple and  lemon transforms which induce translation of the target function, and discussed an example machine geometry from airport baggage screening, first introduced in \cite{WebberQuinto2020II}. We analyzed two lemon transforms, namely $\mathcal{L}_T$ and $\mathcal{L}_0$ (see Corollary \ref{corr1} and \eqref{def:AoLo} respectively), which were shown to satisfy the Bolker condition when the function support was restricted to the upper half of the unit ball. The corresponding apple transforms $\mathcal{A}_T$ and $\mathcal{A}_0$  were shown to violate the Bolker condition. Specifically, there were artifacts induced on the intersections of apples and cylinders with the same axis of revolution. This indicates higher instability in $\mathcal{A}_T$ and $\mathcal{A}_0$ inversion, when compared to $\mathcal{L}_T$ and $\mathcal{L}_0$. Thus, it may be beneficial to use forward scattered photons, which correspond to lemon integrals, if one were to manufacture a CST machine with linear scanning motion (e.g., a scanner with translated sources and detectors, as considered here). The theory is not global however, and does not account for all CST geometries which include linear motion. In further work, we aim to generalize our theory and determine whether such artifacts as discovered here are present for any CST modality with translated sources and detectors.

\section*{Acknowledgments}
This material is based upon work supported by the U.S. Department of Homeland Security, Science and Technology Directorate, Office of University Programs, under Grant Award Number 70RSAT19FR0000155. The views and conclusions contained in this document are those of the authors and should not be interpreted as necessarily representing the official policies, either expressed or implied, of the U.S. Department of Homeland Security.

The second author thanks the U.S. National Science foundation for grant DMS 1712207 and Simons Foundation for grant 708556 that partially supported this research. 

\bibliographystyle{abbrv} 
\bibliography{RefMicGen2}

\end{document}